\documentclass[11pt]{article}
\usepackage{amsmath}
\usepackage{amssymb}
\usepackage{amsfonts}
\usepackage{graphicx}
\usepackage{epsfig}
\usepackage{color}
\usepackage{comment}
\usepackage{epsfig,amsmath,amssymb,latexsym,color,placeins}
\usepackage{algorithmic}
\usepackage[ruled,vlined]{algorithm2e}
\usepackage{graphicx}
\usepackage[colorlinks]{hyperref}
\usepackage{hyperref}
\hypersetup{colorlinks=true, pdfstartview=FitV, linkcolor=black, citecolor=black, plainpages=false, pdfpagelabels=true, urlcolor=black}
\usepackage{subcaption}
\usepackage{xcolor}
\usepackage{epstopdf}
\usepackage[toc,page]{appendix}
\usepackage{siunitx,array,multirow}

\def\bfc{{\bf c}}

\def\bff{{\bf f}}

\def\bfw{{\bf w}}

\def\bfx{{\bf x}}
\def\bfy{{\bf y}}

\def\det{{\rm det}}

\def\MAE{Monge-Amp\'ere equation }

\newtheorem{theorem}{Theorem}
\newtheorem{lemma}{Lemma}

\newtheorem{corollary}{Corollary}

\newtheorem{example}{Example}

\newenvironment{proof}{\begin{trivlist}\item[]{\emph{Proof.}}}
      {\hfill$\Box$\end{trivlist}}

\begin{document}
\title{A Bivariate Spline based Collocation Method for\\Numerical Solution to Optimal Transport Problem}
\author{Ming-Jun Lai\footnote{mjlai@uga.edu, Department of Mathematics, University of Georgia, Athens, GA 30602.
This author is supported
by the Simons Foundation Collaboration Grant \#864439.}
\and Jinsil Lee
\footnote{Department of Mathematics, Ewha Womans University, Seoul, 03760, South Korea, jl74942@ewha.ac.kr}}
\date{}
\maketitle

\centerline{\bf Abstract}
In this paper, we study a spline collocation method for a numerical solution to the optimal transport problem 
 We mainly solve the \MAE with the second boundary condition numerically by  
proposing a center matching algorithm. We prove a
pointwise convergence of our iterative algorithm under the assumption the boundedness of spline iterates. 
We use the \MAE with 
Dirichlet boundary condition and some known solutions to the \MAE with second boundary condition to demonstrate the 
effectiveness of our algorithm. Then  
we use our method to solve some real-life problems. One application problem is to use the optimal transportation 
for the conversion of fisheye view images into standard rectangular images.  
 
\bigskip
\leftline{{\bf Keywords and phrases:} Spline Collocation Method, Monge-Amp\'ere equation, }
\leftline{\hspace{2in} Optimal Transport Problem}

\leftline{{\bf Mathematics Subject Classification:} 35J96, 35J25, 35B65, 65N30, 65K10}

\section{Introduction}
The well-known problem of optimal transportation can be described as follows:
We start with a pile of sand over a domain $V \subset \mathbb{R}^d$ (where $d>1$) and we 
aim to move it to another location $W \subset \mathbb{R}^d$. The volume of the pile will 
be the same before and after the movement, which we can assume to be one for simplicity. 
The shape (or density) of the pile is given by the functions $f$ and $g$ over $V$ and 
$W$, respectively, such that the integrals of $f$ and $g$ over $V$ and $W$ are equal to 
one. Thus,
\begin{equation}
\label{vol}
\int_V f(\bfx) d\bfx =  \int_W g(\bfw) d\bfw =1.
\end{equation}
To move the pile $f$ over $V$ to become pile $g$ on $W$,  we need to find a transform 
$\bfy = T(\bfx) \in W$ for all $\bfx \in V$ that will change the density $f(\bfx)$ at $\bfx \in V$ to the density 
$g(\bfy) = g(T(\bfx))$ at $\bfy \in W$. We can represent these transforms as a vector value function over $V \times W$, 
and the collection of all such mappings that transform $f$ over $V$ to $g$ over $W$ is denoted by $\Pi(V,W)$. For 
any $p \in \Pi(V,W)$, $p$ satisfies $p(A \times W) = f(A)$ and $p(V \times B) = g(B)$ for all open sets $A \subset V$ 
and $B \subset W$.

Moving the sand around requires some effort, which we model by using a cost function $c(\bfx,\bfy)$ defined over $V 
\times W$, where $c(\bfx, \bfy) \geq 0$ for convenience. The goal of the optimal transport problem is to find the 
movement of the pile over $V$ to a new pile over $W$ at a minimal cost. We define the energy or total cost of moving the 
pile $f$ over $V$ to the pile $g$ over $W$ using the transform method $p$ as $E(p) = \int_{V \times W}c(\bfx,\bfy) d 
p(\bfx,\bfy)$ for all $p \in \Pi(V,W)$.
The Kantorovich optimal transportation problem, named after Leonid Kantorovich who won a Nobel prize in Economic Science 
in 1975 for his related work, i.e. linear programming  in 1945, can be stated as the following minimization problem:
\begin{equation}
	\label{Kproblem}
	\min_{p\in \Pi(V,W)} E(p).
\end{equation}

We can reformulate the optimal transportation problem according to  Monge's original optimal transportation problem, 
which seeks to minimize the cost of moving $f$ to $g$ by finding a transform $T \in \mathcal{T}$ that pushes the density 
of $f$ over $V$ forward to $g$ over $W$. This can be represented as the following minimization problem:
\begin{equation}
	\label{Mproblem}
	\min_{T\in {\cal T}}\int_V c(\bfx,T(\bfx)) f(\bfx) d\bfx,
\end{equation}
where ${\cal T}$ is the collection of all push-forwards, denoted by $g = T\# f$. Monge's optimal transportation problem was first formulated by Gaspard Monge in 1781 and gained fame in 1885 due to the prize offered by 
the Academy of Paris (cf. \cite{V03}). It is known that the set $(\bfx, T(\bfx))\in V\times W$ is a cyclic monotonic 
set under the cost function $c(\bfx,\bfy) = \frac{1}{2} \|\bfx - \bfy\|^2$ if $T$ is a minimizer of (\ref{Mproblem}).  
By Rockafellar's theorem, the set is the subdifferential field of a convex function $u$.     
Then the following theorem is known as Brenier's theorem in 1987.

\begin{theorem}[cf.  \cite{V03}]
\label{BrenierThm1}
Suppose that $f$ is a smooth function over $V$ or a positive density function which does not give a mass to small 
sets. There exists a unique push-forward 
$\nabla u$ with a convex function $u$ which is the minimizer of (\ref{Mproblem}) the cost function $c(\bfx,\bfy) = \frac{1}{2} \|\bfx - \bfy\|^2$.  Furthermore, $u$ satisfies a Monge-Amp\'ere equation:
\begin{equation}
\label{MAE}
\det(D^2u(\bfx)) = \frac{f(\bfx)}{g(\nabla u(\bfx))}, \forall \bfx\in V
\end{equation}
with a second boundary condition
\begin{equation}
\label{sbdc}
\nabla u(\bfx)= \bfy \in  W,~ \forall \bfx\in  V.
\end{equation}
\end{theorem}

Solving the Monge-Amp\'ere equation (MAE) poses significant challenges. 
It is a  nonlinear equation, as indicated by the determinant of the Hessian matrix, given by:
$$
\det(D^2 u)=u_{xx} u_{yy} - (u_{xy})^2
$$
for the 2D case, and
$$
\det(D^2 u)=u_{xx} u_{yy} u_{zz} + 2u_{xy} u_{yz} u_{xz}-u_{xx}(u_{yz})^2- u_{yy} (u_{xz})^2-u_{zz} (u_{xy})^2.
$$
Furthermore, the equation is fully nonlinear as (\ref{MAE}) should 
be rewritten as 
\begin{equation}
\label{NLMAE}
g(\nabla u(\bfx)) \det(D^2u(\bfx)) = f(\bfx), \forall \bfx\in V.
\end{equation} 
In addition, $u$ must be a convex function. We have to impose some
sufficient conditions when solving $u$ to ensure the convexity of 
$u$. For example, when both $f$ and $g$ are strictly positive, 
the equation (\ref{MAE}) will make sure that the determinant of 
Hessian matrix is positive. If we require $\Delta u\ge 0$, then $u$ will be convex in the 2D setting.

Finally, there is another difficulty to solve (\ref{MAE}): the boundary condition (\ref{sbdc}) is not easy 
to impose. The second boundary condition $\nabla u(V)= W$ is sometimes also called 
the natural boundary condition 
(to the optimal transport problem), or oblique boundary condition (as 
the obliqueness was studied in \cite{TW09}), or transport boundary condition in the literature.

When using the Dirichlet boundary condition
\begin{equation}
    \label{DBC}
    u(x) = h(x),  x\in {\partial V}
\end{equation}
for a given function $h$ defined on $\partial V$ only, many numerical methods have been proposed and experimented 
in the past two decades. For example, various finite difference methods have been  proposed and studied in, 
e.g. \cite{CS08},  \cite{BFO10, FO11}, \cite{LG21}, etc. 
 In \cite{FO11}, the researchers used a wide stencil finite difference discretization for the Monge-Amp\'ere equation 
and proved that the solution of the scheme converges to the unique viscosity solution of the equation. They used a 
damped Newton's method to get a solution and proved the convergence of Newton's method. 
They presented numerical results 
for not only smooth solutions but also non-differentiable solutions. 
Several mixed finite element methods are studied in \cite{A13}, \cite{AL14}, \cite{FN09a}, 
as the standard finite element method does not apply because of the nonlinearity and convexity constraint.  
Feng and Neilan used the mixed finite element approximations of the viscosity solution in \cite{FN09a}. By the vanishing 
moment method, they can approximate the original equation by the quasilinear fourth-order PDE. 
They derived error 
estimates for the finite element method using a fixed point and linearization technique. 
In \cite{A13}, it is proved that the numerical solution from a mixed finite element method for 
the Monge-Amp\'ere equation converges to the Aleksandrov solution.
Radial basis functions (cf. \cite{BS19}) and tensor product of 
B-spline functions (cf. \cite{SMRFH19}) are used to solve the \MAE numerically. 
Least squares methods were used in \cite{DG03}, \cite{DG04}, and \cite{CGG18}. The \MAE  is solved by minimizing the 
distance between the Jacobi matrix of $m, Dm$ and a real symmetric positive semidefinite matrix $P$ where $m$ satisfies 
the boundary condition $m(\partial X)=\partial Y$ and $P$ satisfying $\text{det} P=\frac{f}{g(m)}$. 
In \cite{LL23}, a 3D spline collocation method was proposed to solve \MAE with Dirichlet boundary condition (\ref{DBC}) 
which is a very effective and efficient way to find numerical solutions.  

In this paper, we are interested in the numerical solution 
of (\ref{MAE}) with the second boundary condition (\ref{sbdc}). 
In \cite{TW09}, the researchers explained the existence and uniqueness of the solution to (\ref{MAE}) with the second 
boundary condition (cf. Theorem 1.2) and estimated the obliqueness of the solution. A
so-called transport boundary condition was introduced. 
In other words, the boundary condition $\nabla u(\partial V) = \partial W$ can replace the original condition 
$\nabla u(V)=W$ in $V$. More precisely, 
 the researchers in \cite{TW09} showed that when both $V$ and $W$ are convex, the transport condition can be enforced by 
ensuring that the boundary points of $V$ are mapped to the boundary points of $W$. 
In a recent paper (cf. \cite{CLW21}),
the researchers improved the regularity results of the solution of 
(\ref{MAE}) with (\ref{sbdc}) by reducing the requirement of the boundary of the underlying domains of interest. They showed that the solution $u$
is in $C^{2,\alpha}$ when $f>0$ is in $C^\alpha(\bar{V})$ over the closure of $V$ under the assumption of that both $V$ 
and $W$ are convex with $C^{1,1}$ boundary (cf. Theorem 1.1 of \cite{CLW21}).   See the classic regularity 
results mentioned therein.  

In \cite{BB00}, the researchers formulated the optimal transport problem as a time-dependent fluid flow problem. 
See also \cite{AHT03} for a similar fluid flow approach based on gradient descent method. The researchers in 
\cite{HRT10} provided another efficient numerical method for the solution of the optimal transport problem. 
In \cite{LR17}, the researchers discretized 
the transportation and used the nice property of cyclic monontoncity to form their computational procedure. 
In \cite{BFO10}, two numerical methods for \MAE with Dirichlet boundary condition and second boundary condition were introduced. 
In particular, their method called subharmonic iteration looks simple and is effective.
\begin{equation}
\label{alg}
\Delta u_{k+1} = \sqrt{ (\Delta u_k)^2 + 2(f/g(\nabla u_k) - \det(D^2 u_k))}
\end{equation}
to solve the \MAE with the Dirichlet boundary condition. 
When solving the \MAE with the second boundary condition 
the  researchers in  \cite{F12} and \cite{BFO14} introduced the projection technique by 
 mapping the $\nabla u(x)$ onto $\partial W$ for all $x\in \partial V$.  
The researchers used finite difference discretization and the system of equations was solved by using Newton's 
iterative method.  Note that the convergence of finite difference methods for the \MAE was established in \cite{CS08}.   
Their computational results demonstrated the effectiveness of the projection method 
although the method is not able to 
find the solution of all optimal transport problem as explained in their paper.  
We note that when $\Omega$ is not convex, the projection $p^k$ is not unique and hence is not well defined. 
Certainly, least squares methods were also used to solve the optimal transport problem, see, e.g. \cite{PBRTIT15A}.


This paper proposes to use bivariate splines for numerical solutions to the optimal transport problem due to the
theoretical study in \cite{CLW21}. That is, when $f>0$ and $f\in C^\alpha(\bar{V})$, 
the solution to the \MAE with boundary condition (\ref{sbdc})
is in $C^{2,\alpha}(\bar{V})$ which allows us to use $C^2$ splines to approximate the solution.  
Furthermore,  our spline collocation method has  
several advantages. Firstly, we can compute the numerical solution 
over an arbitrary polygonal in the 2D setting and polyhedral domain in the 3D setting 
(cf. \cite{LL22} and \cite{LL23}). In 
particular, we are able to solve the optimal transport problem over star-shaped domains as explained in this paper. 
Secondly, the collocation method is much simpler than the weak formulation as seen in \cite{FN09a}, 
\cite{A13}, \cite{A14}, etc. Thirdly, the collocation method can deal with domains with piecewise curve boundaries as explained in \cite{S19}. 
For theoretical properties and numerical implementation of bivariate/trivariate spline functions, 
see \cite{LW04}, \cite{ALW06}, \cite{LS07},  \cite{S15}, \cite{LL22} as well as several dissertations 
(\cite{M19}, \cite{X19} and \cite{L23}) written 
to explain how to implement and how to use multivariate splines for numerical 
solutions of Helmholtz equations, Maxwell equations, and 3D surface 
reconstruction.

We mainly propose a boundary computational method which will be called the center 
matching method.  See Algorithm~\ref{alg2} given in \S 3 with some heuristic reasons.  
To solve the \MAE numerically, we use the same iterative method in (\ref{MAEalg}).
This iterative method has already  been used in several recent works (cf. \cite{A16},\cite{LG21}, \cite{LL23} and etc.).  
It is easy to see that $\Delta u$ is the fixed point of the iterative algorithm as pointed out in \cite{BFO10}. Also, the
iterative is well defined as explained in \cite{A16}, i.e. the term inside the square root is nonnegative 
by using the inequality 
\begin{equation}
\label{Awanou}
\det(D^2 u_k) \le \frac{1}{n^n} (\Delta u_k)^n
\end{equation}
for $n\ge 2$.  
However, the proof of the convergence of our iterative Algorithm~\ref{MAEalg} and Algorithm~\ref{alg2} is still elusive. 
We should provide some  justification for a point-wise convergence under the assumption that the iterates are uniformly 
bounded and an average convergence if the iterates are moderately unbounded.     

\begin{figure}[thpb]
\centering 
\includegraphics[width = 0.75\textwidth]{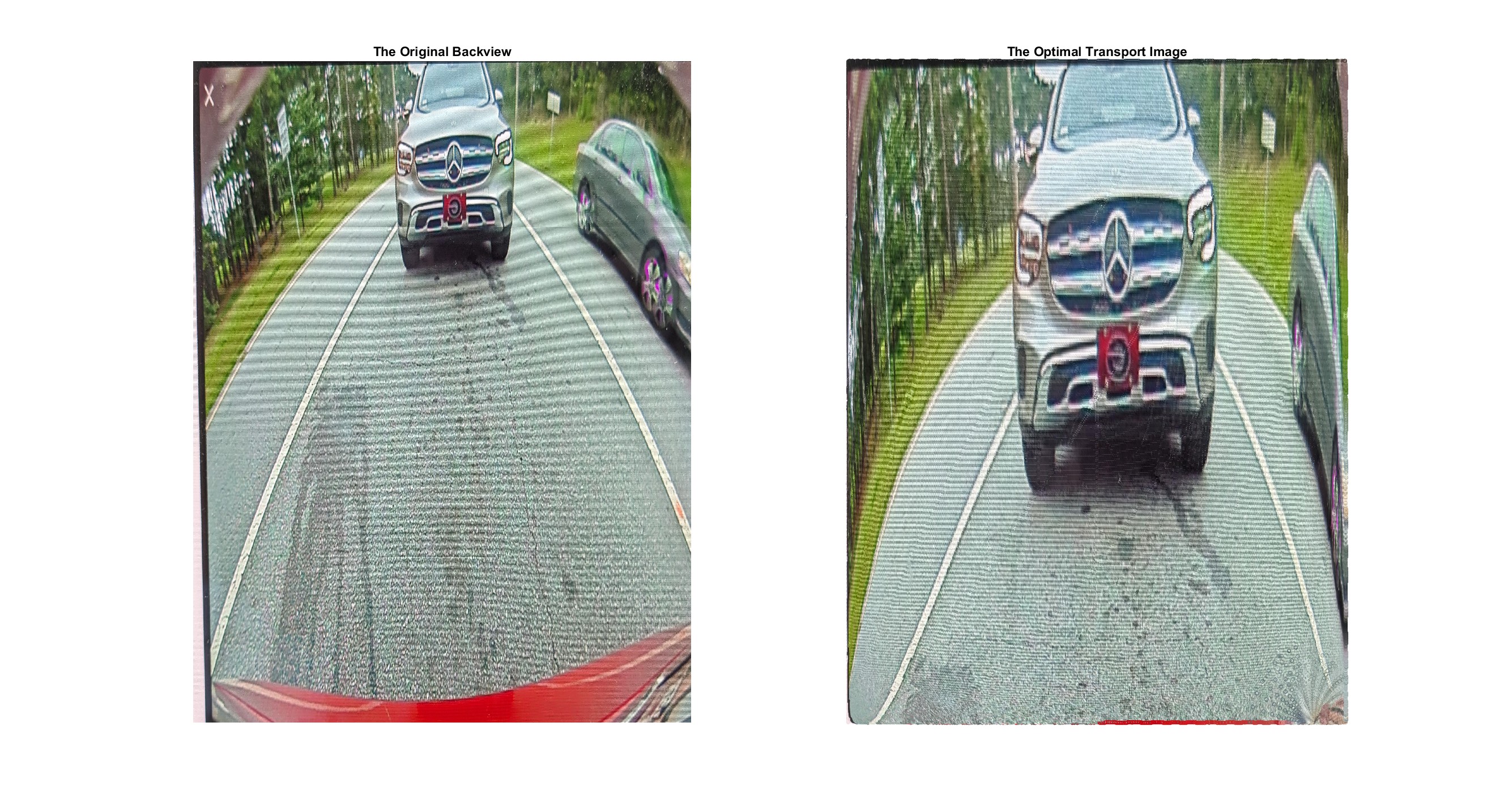} 
\caption{The standard back view of a red car (left) which looks far away from the red car with the backview camera 
	and the optimal transport view (right) by the method in this paper which show that the cars in the back are truly close to the red car. }
\label{backview}
\end{figure}

We first demonstrate that our bivariate spline collocation method can successfully solve 
the \MAE with Dirichlet boundary condition.  Then  we explain our center matching algorithm together with 
the iterative algorithm to compute a spline approximation of 
the solution to the optimal transport problem. 

\begin{figure}[thpb]
\centering 
\includegraphics[width =1\textwidth, height = 0.5\textwidth]{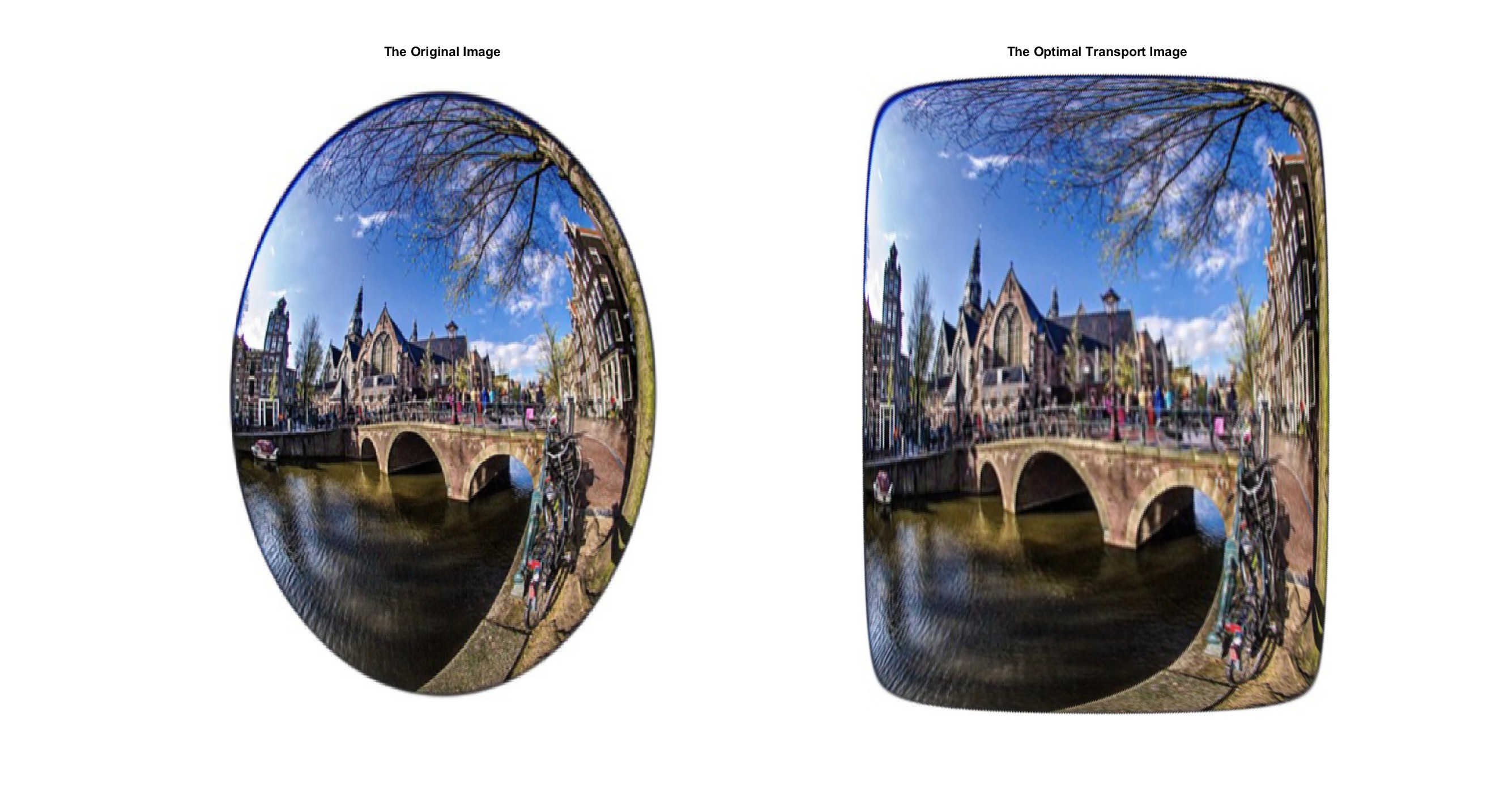} 
\caption{The fisheye view (left) and the optimal transport view (right) by using the method in this
paper. }
\label{fishview}
\end{figure}

Finally we use our computational method to solve
the optimal transport problem which does not have a known solution in general. In particular, we solve 
some real life problems as shown in Figure~\ref{backview} and Figure~\ref{fishview}. Indeed, many cars use back cameras 
to show the back view of the car. The view is not realistic in the sense that the cars in the back look far away 
from the car with back cameras as shown on the left of Figure~\ref{backview}. 
It is interesting and may be necessary to make the view close to the real life situation 
as demonstrated on the right of Figure~\ref{backview}. 

For another example, there are many fisheye images due to the surveillance cameras.  It is interesting and may be necessary to convert these images into a normal standard view.  We use the optimal
transport method to do the conversion of these images as shown in Figure~\ref{fishview}. 

\begin{figure}[thpb]
\centering 
\includegraphics[width =1\textwidth, height = 0.5\textwidth]{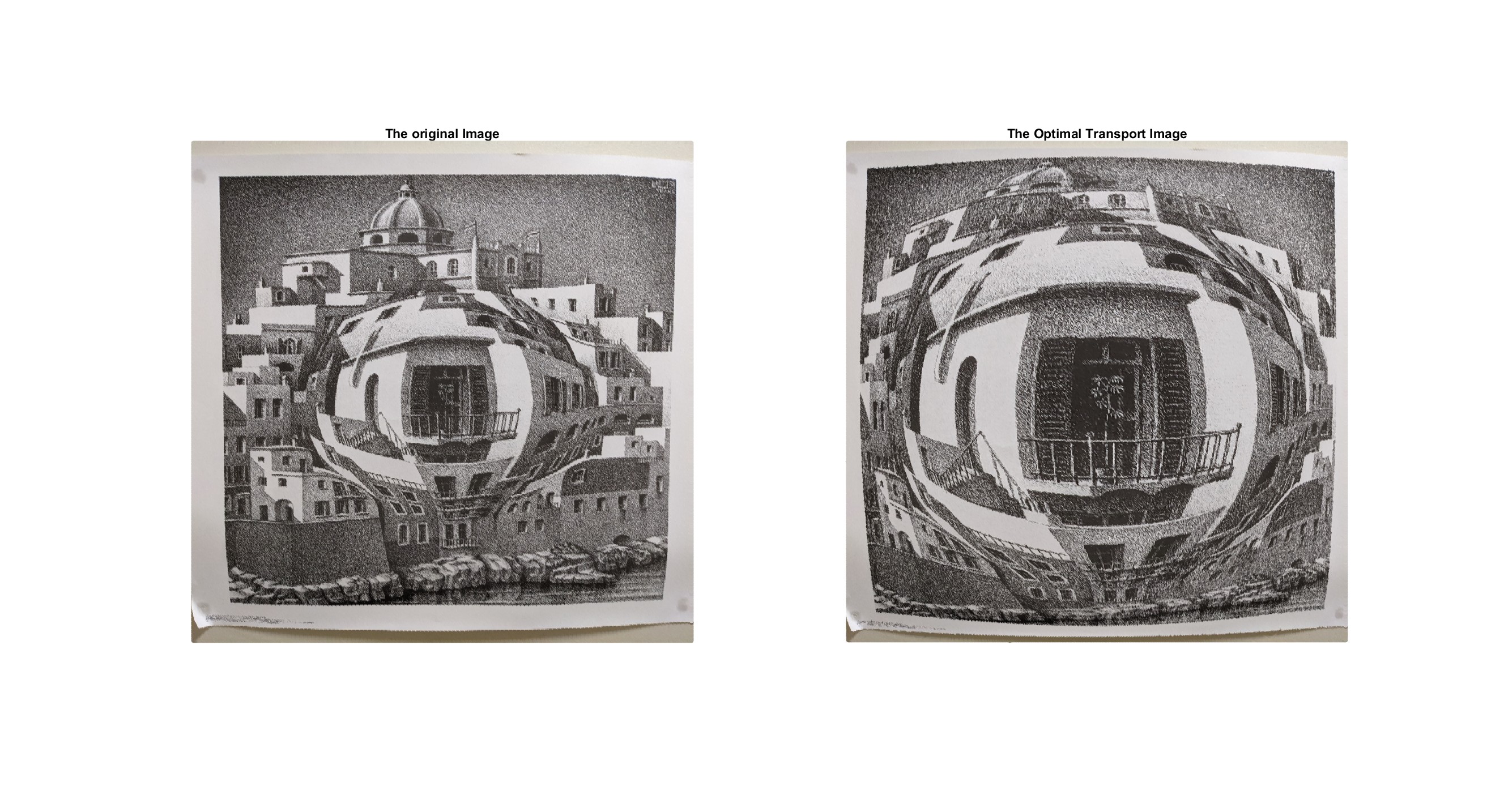} 
\caption{An Escher's view of buildings (left) and the optimal transport view (right) by using the method in this
paper. }
\label{Escherview}
\end{figure}

Motivated by one of Escher's arts (see the left image of Figure~\ref{Escherview}, we use the solution 
to optimal transport problem to generate the image on the right of Figure~\ref{Escherview}. More examples can be found
in the last part of this paper.

We shall demonstrate that our spline solutions $u$ for the optimal transport view in
these figures above are convex,  the gradient map 
$\nabla u$ transforms $V$ into $W$ with $\nabla u(V)= W$  
illustrated by using images above and the RMSE of the \MAE is very small in tne end of this 
paper. By using the uniqueness of the solution in Brenier Theorem above, we conclude 
that the spline solution is a good approximation solution to the optimal transport problem.  
Many numerical examples will be shown to 
demonstrate that our spline method is effective and efficient.

\section{A Spline Collocation Method for Monge-Amp\`ere Equation}
Let us begin the following known result.
\begin{theorem}[Chen, Liu, Wang, 2021\cite{CLW21}] 
\label{Bound}
Assume that $V$ and $W$
are bounded convex domains in $\mathbb{R}^n$  with $C^{1,1}$ boundary and assume that 
$f \in C^\alpha(\bar{V})$ is positive with $0<\alpha<1$. Let $u$ be a convex solution to (\ref{MAE}) and 
(\ref{sbdc}). Then we have the estimate
\begin{equation}
\|u\|_{C^{2,\alpha}(\bar{V})} \le  C,
\end{equation}
where $C$ is a constant depending only $n, \alpha, f, V, W$. Furthermore, if $\alpha=0$, i.e. $f\in C^0(\bar{V})$, then  
\begin{equation}
\|u\|_{W^{2,p}(V)} \le  C
\end{equation}
for all $p\ge 1$, where $C$ is another constant dependent on  $n, p, f, V, W$.
\end{theorem}

The smoothness of the solution of the \MAE  enables us to use $C^2$ smooth spline functions to solve the 
\MAE.    We propose the following computational 
algorithm for the \MAE with Dirichlet boundary condition
called the subharmonic 
iteration method which is a  modified version of the 
original algorithm in \cite{BFO10}.  

\begin{algorithm}[h]
\caption{Subharmonic Iterative Method}
\label{MAEalg}
\begin{enumerate}
    \item 
We start with an initial guess $u^0$ and solve the following Poisson equation first 
\begin{equation}
\label{u0}
\Delta u_1 = 2\sqrt{f/g(\nabla u^0)} \hbox{ in } V\hbox{ and }  u_1 = h, \hbox{ on }\partial V. 
\end{equation}
For $1\le k\le n_1-1$ with $n_1>1$, repeat the following computation until a stopping condition is met:
\begin{equation}
\label{Iterations}
\Delta u_{k+1} = \sqrt{ (\Delta u_k)^2 + 4(f/g(\nabla u^0) - \max\{0,\det(D^2 u_k)\}) }, \hbox{ in } V\hbox{ and }
 u_{k+1}  = h, \hbox{ on }\partial V.
\end{equation}
Stop if $\|\Delta u_{k+1}\|$ is bigger than a fixed amount or if $\|u_{k+1}- u_k\|\le \epsilon$ 
for a given tolerance $\epsilon>0$. 

\item Next for $j= 1$. We use $u^{n_1}$ to find the following new iterations:
\begin{equation}
\label{u1}
\Delta u_{n_1+1} = 2\sqrt{f/g(\nabla u^{n_1})} \hbox{ in } V\hbox{ and } u_{n_1+1} = h, \hbox{ on }\partial V.
\end{equation}
For $1\le k\le n_{j+1}-1$ with $n_{j+1}> 1$, we solve 
\begin{equation}
\label{Iterations1}
\left\{ \begin{array}{cl}
\Delta u_{n_1+k+1} &= \sqrt{ (\Delta u_k)^2 + 4(f/g(\nabla u^{n_1}) - \max\{ 0, \det(D^2 u_{n_1+k})\})}, \hbox{ in } V\cr 
u_{n_1+k+1}  &= h, \hbox{ on }\partial V. 
\end{array}
\right.
\end{equation}
Stop if $\|\Delta u_{n_1+k+1}\|$ is bigger than a fixed mount or if $\|u_{n_1+k+1}- u_{n_1+k}\|\le \epsilon$ 
for a given tolerance $\epsilon>0$. 

\item Repeat the above process for $j=2, \cdots, N$. 
\end{enumerate}
\end{algorithm}


As explained in Appendix,  
our collocation method is based on the discretization of (\ref{Iterations}) and (\ref{Iterations1}) by using bivariate splines and then finds 
${\bf c}^*$ which solves the following constrained minimization: 
 \begin{align}
\label{min1ma}
\min_{\bf c} J(c)=\frac{1}{2}(\alpha\|B{\bf c}^{k+1} -  G \|^2+ \beta \|H{\bf c}^{k+1} \|^2)  \text{  subject to } 
K{\bf c}^{k+1}= {\bf f}^{k},
\end{align}
where $B, G$ are from the boundary condition, $H$ is a matrix from the smoothness conditions of spline functions (cf. \cite{LS07})  and $\bff^k$ is the right-hand side
of (\ref{Iterations}) dependent $u_k$ and $K$ is the matrix associated with the discretization of the Laplace operator.

 According to \cite{A14}, it is known
that if  $\text{det} (D^2 u^*)=f>0$ and $u^*$ is convex, then there exists a neighborhood of $u^*$,
say $N(u^*)$ 
such that the eigenvalues of the Hessian matrix 
$(D^2 u(x))$ are nearby the 
eigenvalues of the Hessian matrix $D^2 u^*(x), \forall x\in V$ for 
all $u\in N(u^*)$. 
In other words, 
there exist 
two positive constants $m, M$ dependent on the eigenvalues of $D^2u(x)$ such that 
$$0<m \le \lambda_2\le \lambda_1 \le M, $$
where $\lambda_1,\lambda_2$ are the eigenvalues of 
$ (D^2 u(x)), \forall x\in V$ for any $u\in N(u^*)$. 
Next, the following result is also known (cf. \cite{A15}). For clarity, we provide proof below. 
\begin{lemma}
	\label{eigen1}
	Suppose that the convex solution $u^*\in W^{2,\infty}$ satisfies $\text{det} (D^2 u^*)
	=f/g(\nabla u^*)>0$. There exists a $\delta>0$ such that  for any $u$ which is close 
	enough to the exact solution $u^*$ in the sense that $|u-u^*|_{2,\infty}\le \delta$, we have, for any $a<4$,  
	\begin{equation*}
		\text{det} (D^2 u)\le \frac{1}{4}(\Delta u)^2< \frac{1}{a}(\Delta u)^2.	\end{equation*} 
\end{lemma}
\begin{proof}
	It is clear that
	$$\lambda_1 \lambda_2 =\text{det} (D^2 u)\le \frac{1}{4}(\Delta u)^2=\frac{1}{4}(\lambda_1+\lambda_2)^2$$
	\end{proof}
 
Using this Lemma \ref{eigen1}, we can conclude that 
\begin{eqnarray*}
\Delta u_{k+1}&=&\sqrt{(\Delta u_k)^2+4(f/g(\nabla u^j) -\det D^2 u_k)}\cr
&=& \sqrt{(\Delta u_k)^2-4\det D^2 u_k+4f/g(\nabla u^j)}\geq 2\sqrt{f_{min}/g_{max}}>0,
\end{eqnarray*}
where $f_{min}=\min\{ f(\bfx), \bfx\in V\}$ and $g_{max} = \max\{g(\bfy): \bfy\in W\}$.  That is, 
$\Delta u_{k+1} >0, k\ge 1$. 
If we also have $\det(D^2 u_k)>0$, we know $u_k$ is convex 
as both eigenvalues of the hessian $D^2u_k$ have the same signs and 
the addition of the two eigenvalues is positive.  

We shall make the following assumption for our convergence analysis. 

\noindent
{\bf Assumption 1.} Both $f(x,y)\ge f_0>0$  over $(x,y)\in V$ and $g_{max}\ge g(x,y)\ge g_0 >0$ over $(x,y)\in W$. 

\noindent
{\bf Assumption 2.}   $\det(D^2 u_k)>0$ for all $k\ge 1$. We shall use $\det(D^2 u_k)$ instead of 
$\max\{0, \det(D^2 u_k)\}$ for simplicity in the discussion later. 

\noindent
{\bf Assumption 3.}  $\Delta u_k$, $k\ge 1$ have an upper bound over $\Omega$. 

Note that we shall show that $\Delta u_k\ge 0$. Together with Assumption 2, we see that $u_k$'s are convex 
in the 2D setting.   
Assumption 3 is needed because of Theorem~\ref{Bound}. 
Also, since we want to have $(\frac{f}{ g(\nabla u^j)}-
\text{det} (D^2 u_{n_j+k}))$ converge to zero for some  $j\ge 1$, 
$\det(D^2 u_{n_j+k})$ must be bounded.  By Lemma~\ref{eigen1}, Assumption 3 implies the boundedness of 
$\det(D^2 u_{n_j+k})$ in the 2D setting.  
Theorem~\ref{mainresult} below says that Assumption 3 is necessary to ensure the point-wise convergence 
of Algorithm~\ref{MAEalg}.
If they are unbounded, we simply choose another initial guess to restart Algorithm~\ref{MAEalg}. 

With these three assumptions, we derive some properties for convergence of our main algorithm~\ref{MAEalg}. 
We first use the subharmonic iterative step \eqref{Iterations} to get 
		\begin{align*}
	(\Delta u_{n_j+k+1})^2  
 &=(\Delta u_{n_j+k})^2
	+4(\frac{f}{ g(\nabla u^j)}-\text{det} (D^2 u_{n_j+k}))\cr
	&=(\Delta u_{n_j+k-1})^2+4(\frac{f}{ g(\nabla u^j)}-\text{det} (D^2 u_{n_j+k-1}))
+4(\frac{f}{ g(\nabla u^j)}-\text{det} (D^2 u_{n_j+k}))\\
	&= \cdots \cr 
	&=(\Delta u_{n_j})^2+4\sum_{l=1}^k\Big{[}\frac{f}{ g(\nabla u^j)}-\text{det} (D^2 u_{n_j+l})\Big{]}
\end{align*}
for all $k\le n_{j+1}$, By Lemma~\ref{eigen1} and the iteration restarting step (\ref{u1}), we have 
$$
4\det(D^2 u_{n_j+k+1}) \le 4\frac{f}{ g(\nabla u^j)} +4\sum_{l=1}^k\Big{[}\frac{f}{ g(\nabla u^j)}-\text{det} (D^2 u_{n_j+l})\Big{]}.
$$
In other words, the following nonnegativity property holds.
\begin{equation}
\label{nonnegativity0}
4\sum_{l=1}^{n_j}\Big{[}\frac{f}{ g(\nabla u^j)}-\text{det} (D^2 u_{n_j+l})\Big{]} 
+ 4\frac{f}{ g(\nabla u^j)}- 4\det(D^2 u_{n_j+k+1}) \ge 0.    
\end{equation}
Then, we have
\begin{equation}
\label{nonnegativity}
\sum_{l=1}^{n_{j+1}}\Big{[}\frac{f}{ g(\nabla u^j)}-\text{det} (D^2 u_{n_j+l})\Big{]}  \ge 0. 
\end{equation}
Similarly, we have
\begin{equation}
\label{nonnegativityall}
\sum_{j=0}^{N}\sum_{l=1}^{n_{j+1}}\Big{[}\frac{f}{ g(\nabla u^j)}-\text{det} (D^2 u_{n_j+l})\Big{]}  \ge 0
\end{equation}
for any integer $N>0$.
The discussion above leads to the following 
\begin{theorem}
\label{mainresult}
	Fix a spline space $\mathcal{S}^r_D(\triangle)$ with $\triangle$ being a triangulation of the domain
	$\Omega$ with smoothness $r\ge 2$ and degree $D\ge 3r+2$. 
	Let  $u_k\in \mathcal{S}^r_D(\triangle), k\ge 1$ be the sequence from Algorithm~\ref{MAEalg}. 
	Then we have the nonnegativity properties \eqref{nonnegativity0} and  \eqref{nonnegativity}. 
 Furthermore, under Assumption 3, we have 
	\begin{equation}
		\label{mjlai03092022}
\sum_{k=1}^{n_{j+1}}(\frac{f(\bfx)}{g(\nabla u^j)} -\det (D^2 u_{n_j+k})(\bfx)) 
\to 0  
	\end{equation}
	when $j  \to \infty$ for all $\bfx\in \Omega$.
\end{theorem} 
\begin{proof}
Since we have for any given $j$
\begin{align*}
 (\Delta u_{n_j+n_{j+1}})^2-(\Delta u_{n_j})^2=4\sum_{l=1}^{n_{j+1}}\Big{[}\frac{f}{ g(\nabla u^{j})}- 
\text{det} (D^2 u_{n_j+l})\Big{]},
\end{align*}
we sum the above equalities for $j=0, \cdots, N$ to have  
$$
(\Delta u_{\sum_{j=0}^Nn_j +1})^2-(\Delta u_{0})^2=4\sum_{j=0}^N \sum_{l=1}^{n_{j+1}}\Big{[}\frac{f}{ g(\nabla u^{j})}- 
\text{det} (D^2 u_{n_j+l})\Big{]},
$$ 
Hence, by  Assumption 3 we conclude (\ref{mjlai03092022}) since the boundedness of the summation on the right-hand side
of the above equation with nonnegative 
terms to finish the proof of Theorem~\ref{mainresult}.
\end{proof}

Another simple conclusion can be made based on the iterations in Algorithm~\ref{MAEalg}. 
\begin{lemma}
\label{mjlai09182023}
Under Assumption 2, we have $\|\Delta u_k\|_\infty\le (k+1)2\| \sqrt{f/g_{min}}\|_{\infty}$ for all $k\ge 0$, 
where $g_{min}=\min\{g(\bfy), \bfy\in V\}$.   
\end{lemma}
\begin{proof}
We use an induction to establish the proof. 
Clearly, for $k=1$, we have $\|\Delta u_1\|_\infty \le \|2\sqrt{f/g_{min}}\|_\infty$.
For a general $k>1$, we have
$$
( \Delta u_{k+1})^2 \le (\Delta u_k)^2 + 4(f/g(\nabla u^0))\le (k+1) 4 \|f/g_{min}\|_{\infty} ^2 + 4(f/g_{min}) 
$$ 
by using the induction hypothesis.  It follows $\|\Delta u_{k+1}\|\le (k+2)2\| \sqrt{f/g_{min}}\|_{\infty}$. 
\end{proof}

That is, the growth of $\Delta u_k$ is at most about $O(n)$. 
If the iterative solution $\Delta u_{n_1+\cdots+n_j+k}$'s are moderately unbounded in the 
sense that $\|\Delta u_{n}\|=o(n)$, then we have the following 

\begin{theorem}
\label{mainresult2}
Suppose $\Omega \subseteq \mathbb{R}^2$ is a bounded  domain. Suppose that  
$\|\Delta u_k\|_{L^2(\Omega)} = o(k)$ for all $k>1$. Then 
$(f/g(\nabla u^j)-\text{det} (D^2 u_{n_j+l})\to 0$ in average in the sense
\begin{equation}
\label{averageconvergence}
\frac{1}{N} \sum_{j=0}^{N}\sum_{l=0}^{k+1}\Big{[}\frac{f}{ g(\nabla u^j)}-\text{det} (D^2 u_{n_j+l})\Big{]}  \to  0
\end{equation}
when $N\to \infty$. 
\end{theorem}

Finally we  establish the following
\begin{theorem}
\label{mainresult2}
Consider the case $g\equiv 1$ over $W$.  
Under assumptions 1--3, there exists a real number $\rho<1$ such that 
the iterative solutions $u_1, \cdots, u_{k+1}$ from (\ref{Iterations}) satisfy 
\begin{equation}
\label{conv}
\|u - u_{k+1}\|_{H^2} \le (\frac{\rho  \sqrt{2}}{A})^k \|u - u_1\|_{H^2}.
\end{equation}
where $A>0$ is a positive constant in Lemma~\ref{normequiv}.  
Furthermore, suppose that $\frac{\rho  \sqrt{2}}{A}\le \rho_0 <1$. Then
the iterative solutions $u_k, k\ge 1$ converge to the exact solution $u$.
\end{theorem}

\begin{proof}
 Let us first rewrite the right-hand side 
$$ \sqrt{(\Delta u_k)^2+4(f/g(\nabla u^j)-\det D^2 u_k)} = \sqrt{ (D_{xx} u_k + D_{yy} u_k)^2+ 4f/g(\nabla u^j)-4 D_{xx} u_k D_{yy} u_k+4(D_{xy} u_k)^2}.$$  As the solution
$u$ is a fixed point of the iteration, we have
\begin{eqnarray*}
&&\Delta u - \Delta u_{k+1} \cr 
&=& \sqrt{ (u_{xx} -u_{yy})^2+ 4(u_{xy})^2+ 4f/g(\nabla u^j)}
- \sqrt{((u_k)_{xx} -(u_k)_{yy})^2+ 4((u_k)_{xy})^2+ 4f/g(\nabla u^j)}\cr 
&=& \frac{(u_{xx} -u_{yy})^2+ 4(u_{xy})^2 -((u_k)_{xx} -(u_k)_{yy})^2- 4((u_k)_{xy})^2}
{\sqrt{ (u_{xx} -u_{yy})^2+ 4(u_{xy})^2+ 4f/g(\nabla u^j)}
+ \sqrt{((u_k)_{xx} -(u_k)_{yy})^2+ 4((u_k)_{xy})^2+ 4f/g(\nabla u^j)}}.
\end{eqnarray*}
Thus, for each point on $\Omega$, we use Cauchy-Schwarz inequality first and then Minkowski inequality to have
\begin{eqnarray*}
&& |(a-b)^2+(2c)^2 -(a_1-b_1)^2-(2c_1)^2| \cr
&= &  |\Big{(}a-b-(a_1-b_1)\Big{)}\Big{(}a-b+(a_1-b_1)\Big{)}+(2c-2c_1)(2c+2c_1)| \cr
&\le &  \sqrt{\Big{(}a-b-(a_1-b_1)\Big{)}^2+(2c-2c_1)^2}\sqrt{\Big{(}a-b+(a_1-b_1)\Big{)}^2+(2c+2c_1)^2}\cr
&\le &  \sqrt{\Big{(}a-b-(a_1-b_1)\Big{)}^2+(2c-2c_1)^2}\sqrt{\Big{(}a-b+(a_1-b_1)\Big{)}^2+(2c+2c_1)^2}\cr
&\le &  \sqrt{(a-b-(a_1-b_1))^2+(2c-2c_1)^2}\Big{[}\sqrt{(a-b)^2+(2c)^2}+\sqrt{(a_1-b_1)^2+(2c_1)^2}\Big{]}\cr
&= & \sqrt{(a-a_1-(b-b_1))^2+4(c-c_1)^2}\Big{[}\sqrt{(a-b)^2+4c^2}+\sqrt{(a_1-b_1)^2+4c_1^2}\Big{]}
\end{eqnarray*}
By letting $a=u_{xx}, b=u_{yy}, c=u_{xy},a_1=(u_k)_{xx}, b_1=(u_k)_{yy}, c_1=(u_k)_{xy}$ and defining
\begin{equation}
\label{keyconst}
\rho(k)=  \frac{ \sqrt{ (u_{xx} -u_{yy})^2+ 4(u_{xy})^2}
+ \sqrt{((u_k)_{xx} -(u_k)_{yy})^2+ 4((u_k)_{xy})^2}} 
{ \sqrt{ (u_{xx} -u_{yy})^2+ 4(u_{xy})^2+ 4f/g(\nabla u^j)}
+ \sqrt{((u_k)_{xx} -(u_k)_{yy})^2+ 4((u_k)_{xy})^2+ 4f/g(\nabla u^j)}}
\end{equation}
which is strictly less than $1$ 
for each $k$ with $f/g(\nabla u^j)\ge f_0/g_m>0$. Assumption 3, i.e. the boundedness of $\Delta u_k$ implies
the boundedness $\Delta (u- u_k)$ in $L_2(V)$ norm.  By the norm equivalence in Lemma~\ref{normequiv}, all $(u_k)_{xx}, 
(u_k)_{xy}, (u_k)_{yy}$ are bounded in $L_2(V)$. Since they are spline functions, they are 
bounded in $L_\infty$ norm. It follows from \eqref{keyconst} that there exists a real number $\rho<1$ such that 
$\rho(k) \le \rho$.  

We summarize the above inequalities together to get 
$$
\|\Delta u - \Delta u_{k+1}\|_{L^2(\Omega)} 
\le \|\rho\|_{L^\infty(\Omega)} \|[(D_{xx} (u-u_k)-D_{yy} (u-u_k)))^2  + 4(D_{xy} (u-u_k))^2]^{1/2}\|_{L^2(\Omega)}.
$$
In terms of our new norm $\|\cdot \|_L$ and standard $H^2$ norm, we have
\begin{equation}
\|u - u_{k+1}\|_L \le \rho \sqrt{2}\|u - u_k\|_{H^2}.
\end{equation}
Now we can use the norm equivalence, c.f. Lemma~\ref{normequiv} again to have 
\begin{equation}
\|u - u_{k+1}\|_{H^2} \le \frac{\rho  \sqrt{2}}{A} \|u - u_k\|_{H^2}.
\end{equation}
Therefore, we finished a proof of the first part of Theorem~\ref{mainresult2}. The second part follows easily. 
\end{proof}

\section{Computational Results}
This section is divided into three parts. 
The first part showcases the accuracy and performance of our numerical method in solving the \MAE with Dirichlet 
boundary conditions. We provide several examples to demonstrate the effectiveness of our approach, using spline 
functions of degree $D=8$ and smoothness $r=2$. Table \ref{MAE:Direx1} shows the RMSE values for the numerical solutions 
$u^{s1}$ of the \MAE on different domains. We observe that our method performs well on all the domains considered.

The second part of this section explains our center matching method for the \MAE with second boundary conditions for 
solving the optimal transport problem. We present several numerical solutions by deforming density functions over convex 
or nonconvex domains. We shall demonstrate that the Brenier potential $u$ is convex, 
the root mean squares error of the \MAE is very small and the gradient map $\nabla u$ is onto 
which imply that the solution obtained is indeed the solution to the optimal transport problem. 

The third part of this section is to apply our Algorithm~\ref{MAEalg} together with the center matching technique for
approximating the solution of some optimal transport problem.

\subsection{Numerical Solution to the \MAE  with Dirichlet Boundary Condition}
We first demonstrate Algorithm~\ref{MAEalg} for a smooth 
solution to show that our spline collocation method works very well.

\begin{example}
In this example, we assess the accuracy of our method in approximating the solution $u^{s1}=e^{\frac{(x^2+y^2)}{2}}$ of 
the Monge Amp\'ere equation with source term $f_1(x,y)=(1+x^2+y^2)e^{x^2+y^2}$ as proposed in \cite{BHP15, BFO10}. 
We use various domains of interest as shown in Figure~\ref{fig:MAEdir} with triangulations  
to test the accuracy of our spline collocation method. We use spline functions of degree $D=8$ and smoothness $r=2$ 
over the triangulations as shown in Figure~\ref{fig:MAEdir} with a total of 20 iterations.	
The RMSE values for these solutions again the exact solution $u^{s1}$ are shown in Table~\ref{MAE:Direx1}. 
The RMSE values are computed based on $201\times 201$ equally-spaced pointed over the bounding box of each 
domain. In addition, numbers of vertices and triangles and computational times in seconds  
are also presented in Table~\ref{MAE:Direx1}. From the table, we can see that our method is able to approximate the testing function $u^{s1}$ very well on all the domains considered. 
In addtion, we have experimented other testing functions, e.g. $z=x\sin(x) + y\sin(y)$ (cf. \cite{FN09a}). 
 Similar  approximation results were produced and the details are omitted here.  

\begin{figure}[htpb]
\centering 
\includegraphics[width=0.2\linewidth]{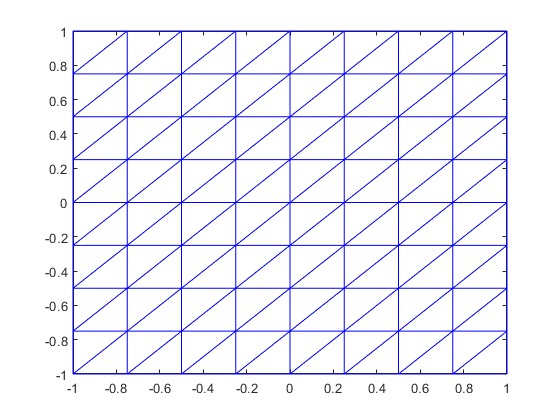}
\includegraphics[width=0.2\linewidth]{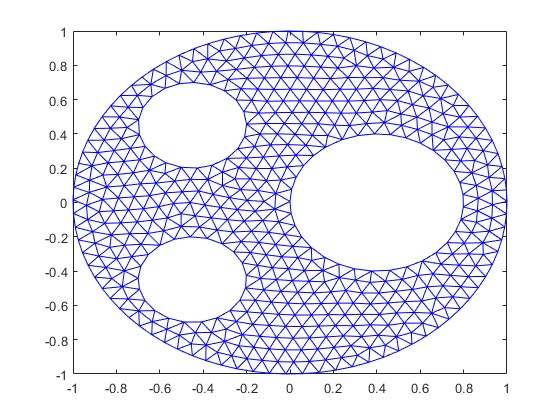}
\includegraphics[width=0.2\linewidth]{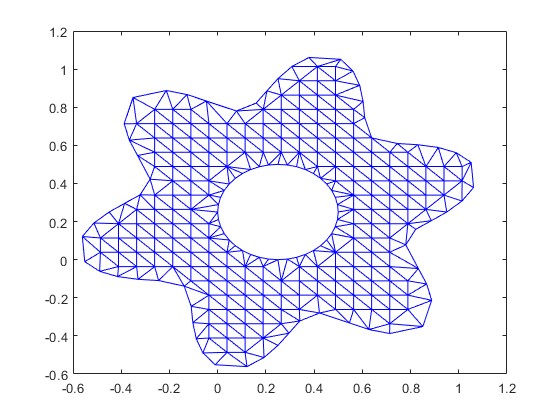}
\includegraphics[width=0.2\linewidth]{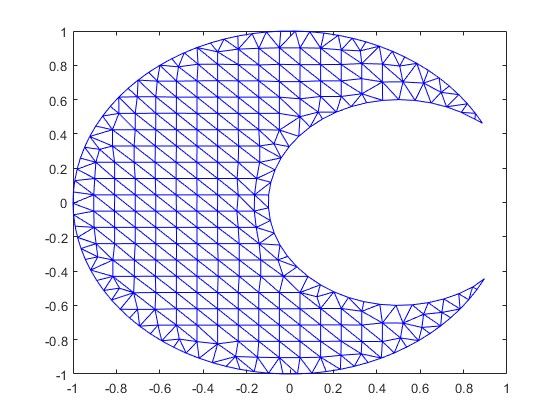}
\caption{2D domains used for solving the Monge Amp\'ere equation with Dirichlet boundary condition
(a) Domain $[-1,1]^2$; (b) Circle with 3 holes; (c) Flower with a hole; (d) Moon }
		\label{fig:MAEdir}
	\end{figure}

\begin{table}[htpb]
	\centering
	\renewcommand{\arraystretch}{1.2} 
	\begin{tabular}{ c c c c c}
		\hline
		Domain & $N_V$ & $N_T$ & CPU time (s) & RMSE \\
		\hline
		$[-1,1]^2$ & 81 & 128 & 1.0937 & 2.6440e-10 \\
		Circle with 3 holes & 525 & 895 & 3.9794 & 1.0894e-10 \\
		Flower with a hole & 297 & 494 & 2.1758 & 3.9749e-10 \\
		Moon & 325 & 531 & 4.8399 & 4.6869e-10 \\
		\hline
	\end{tabular}
	\caption{Numerical results for the Monge Amp\'ere equation with Dirichlet boundary condition on various domains, including the number of vertices($N_V$), number of triangles($N_T$), CPU time (s), and RMSE for spline functions of degree $D=8$ and smoothness $r=2$}
	\label{MAE:Direx1}
\end{table}
 From Table~\ref{MAE:Direx1}, we observe that our method performs well in this domain in Figure \ref{fig:MAEdir}.
\end{example}

We next test our method for the solution which is not very smooth 
and demonstrate that our method still works well.  
\begin{example}\label{ex:3}
We test our method for a solution that is not very smooth and demonstrate its effectiveness. Consider the domain $\Omega:=[0,1]^2$ and the point $x_0:=(0.5,0.5)$.
Consider a solution $u\in C^1(\Omega)$ which is defined by 
$$u^{s3}(x)=\frac{1}{2}(\max\{0, \|x-x_0\|_2 -0.2\})^2 \hbox{ and } f(x)=\max\{0, 1-\frac{0.2}{\|x-x_0\|_2}\}$$
which is a testing function suggested in \cite{BHP15}. The solution is in $C^2(\Omega)$ but has a singularity in its gradient near $x_0$. 
Table \ref{MAE:Direx2} shows the numerical results for the Monge Amp\'ere equation with the above solution, obtained using spline functions of degree $D=8$ and smoothness $r=2$. The table includes the CPU time and RMSE. 
From Table \ref{MAE:Direx2}, we can see that our method performs well. Overall, the results demonstrate the effectiveness of our method for solutions with singularities in their gradients. Further numerical results can be found in \cite{L23} and \cite{LL23}.
\begin{table}[htpb]
	\centering
	\begin{tabular}{ c c c }
		\hline
		\multicolumn{1}{c}{Degree} &CPU time (s)&RMSE \\
		\hline
		8&1.41&3.54e-05 \cr
		\hline
	\end{tabular}
	\caption{Numerical results for approximating the solution to the Monge Amp\'ere equation with singularity, using spline functions of degrees $D=8$ and smoothness $r=2$ over the domain $\Omega=[0,1]^2$}
	\label{MAE:Direx2}
\end{table}
\end{example}

\subsection{Numerical Solution to the \MAE with Transport Boundary Condition}
In this subsection, we propose a computational algorithm, 
Algorithm~\ref{alg2} called a center matching method for solving the \MAE with a transport boundary condition numerically.  

\subsubsection{Our Computational Algorithm}
To motivate our center matching method, let us first look at the solution of
the OTP in simple case with 
$f\equiv 1$.  It is known that when $W= V+z_0$ for a vector $z_0$ and $f\equiv 1$, 
the solution to the OTP is simply the linear translation from $V$ to $V+z_0=W$. That is, $u= \|x+z_0\|^2/2$, 
i.e. $\nabla u=x+z_0$. Indeed, the solution $u$ is convex, $\nabla u$ maps $V$ to $W$ and $u$ satisfies the
\MAE.  For another example, 
if $W=\alpha V$ for $\alpha>0$, then the solution to the OTP is 
$u= \alpha \|x\|^2/2$, i.e., $\nabla u(x) = \alpha x$ if $f=\alpha^2$. Thus, $u$ is convex, and 
$\nabla u$ maps $V$ to $V=\alpha V$, and $u$ satisfies the \MAE with $f'=\alpha^2$, not $f=1$. Note that
the OTP is to minimize $\int_V \|x- y\|^2 f(x)dx = \frac{1}{\alpha^2}\int_V \|x-y\|^2 f'(x)dx$. The minimizers are the same for both minimizing functionals $\int_V \|x- y\|^2 f(x)dx$ and $\int_V \|x- y\|^2 f'(x)dx$. 
As $\nabla u=(\alpha x, \alpha y)$, the geometrical meaning of $\nabla u$ is to push any $(x,y)\in V$ to $(\alpha x, \alpha y)\in W$. In particular, for $(x,y)\in 
\in \partial V$, we have $(\alpha x, \alpha y)\in \partial W$.  This idea motivates us to define our center matching method.  

To deal with the general domain, let us recall the definition of star-shaped domains.
That is, a domain $V$ is star-shaped if there is a point $x_0\in V$ such that 
any ray starting from $x_0$ intersects the 
boundary of $V$ once and only once. Such a point $x_0$ is called a center of $V$. 
Note that any convex domain is a star-shaped domain and any interior point $x_0\in V^\circ$ is a center.   

Let $x_0\in V, y_0\in W$  be  centers of domain $V$ and domain $W$. We may first 
assume that $x_0=y_0$ for simplicity. 
Let us now describe our algorithm as follows. 
For any point $v_\theta \in \partial V$, let $r_\theta$ 
be the ray starting from $x_0$ and passing $v_\theta$. Note that 
$r_\theta$ intersects the boundary $\partial W$ only once and $v_\theta$ is uniquely
determined for any $\theta\in [0, 2\pi)$.  
Let $w_\theta$ be the intersection point of the ray $r_\theta$ and 
$\partial W$. This forms the boundary condition in
(\ref{newboundarycondition}). We now return to the 
general setting. We may first do 
a linear transportation, moving $V$ to $V'=V+y_0-x_0$ which has the cost $\|x_0- y_0\|^2 A(V)$, 
where $A(V)$ is the area of $V$.  Then $V'$ and $W$ have the 
same center and we apply our  Algorithm~\ref{alg2}.
Of course, if $y_0\not=x_0$, we can 
also move $W$ to $W'=W-y_0+x_0$ if
the linear transportation cost 
$\|y_0- x_0\|^2 A(W)$ is smaller than
$\|y_0-x_0\|^2 A(V)$, Now a shift of $W$ 
will have the same center with $V$. Let us assume that W is moved to 
$W'$. For $v_\theta\in 
\partial V$, let $w_\theta= W'_\theta - y_0+ x_0$ in 
(\ref{newboundarycondition}).

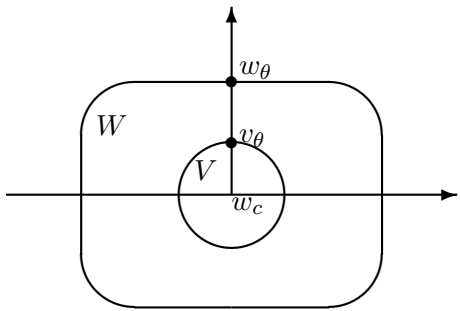
\begin{figure}[htpb]
\setlength{\unitlength}{1cm}
\thicklines
\begin{picture}(20,5)(-1,0)
	\put(6,2){$w_c$}
	\put(4.2,3){$W$}
	\put(5.5, 2.4){$V$}
	\put(6,2.2){\circle{2}}
	\put(6,2.2){\oval(4,3)[r]}
	\put(6,2.2){\oval(4,3)[l]}
	\put(6,2.2){\vector(0,1){2.5}}
	\put(6.1,3.8){$w_\theta$} 
	\put(5.9,3.6){$\bullet$}
	\put(5.9,2.8){$\bullet$}
	\put(6.1,2.9){$v_\theta$}
	\put(3,2.2){\vector(1,0){6}}
\end{picture}
\caption{Distance $\|v_\theta - w_\theta\|^2$ is the smallest for any other $w\in \partial W$, 
where $\theta=90^\circ$ is shown, $V$ is the circular domain and $W$ is a rectangular domain with four round corners\label{ex1}}
\end{figure}

For convenience, in the remaining of this subsection, 
we assume that both $V$ and $W$ are a star-shaped domain and  
$V$ and $W$ share the same center $x_0=y_0$. 
Next, let us give our heuristic reason that  $w_\theta$  chosen in such a way 
can minimize the  transportation cost:   
\begin{equation}
    \int_V | \bfx - T(\bfx)|^2 f(\bfx) d\bfx.
\end{equation}
When  $x_0=y_0$, We need to make the distance $|\bfx- T(\bfx)|^2$ inside of the integral above 
as small as possible. The integrand can be rewritten as $|\bfx- T(\bfx)|^2=\|(\bfx-x_0)- (T(\bfx)-y_0)\|^2$. The difference between two vectors $(\bfx-x_0)$ and $(T(\bfx)-y_0)$ will be small when
the angle between the two vectors is $0$. That is, for three points $x_0, \bfx, 
T(\bfx)$, the distance $|\bfx- T(\bfx)|^2$ will be small when all three points
are in the same line as shown in Figure~\ref{ex1}. That is why we use the same
ray to find the intersection of $v_\theta\in \partial V$ and $w_\theta \in \partial 
W$ for $\theta \in [0, 2\pi)$.

\begin{algorithm}[t] 
\caption{The Center Matching Method\label{alg2}}
\begin{enumerate}
\item Assume that $V$ and $W$  have the same center point $w_c\in W$. 
\item For given $\nabla u^k$ and $v_\theta\in \partial V$, determine the intersection point 
$w_\theta \in \partial W$ by the ray passing $w_c$ and $\nabla u(v_\theta)$ for all $\theta\in 
[0, 2\pi)$.   
\item 
We solve the iterative function $u^{k+1}$ that satisfies the following equations:
\begin{eqnarray}
\text{det} (D^2u^{k+1}(x))&=&\dfrac{f(x)}{g(\nabla u^{k}(x))}, \quad x\in V 
\label{iterativeMAE}\\
\nabla u^{k+1}(v_\theta) \cdot n_x &=&w_\theta \cdot n_x, 
\forall v_\theta \in \partial V \label{newboundarycondition}\\
\int_V u^{k+1} dx&=&0 \label{normalization}
\end{eqnarray}
by using Algorithm~\ref{MAEalg} with the Dirichlet boundary condition replaced by the Neumann boundary conditions above. 
\end{enumerate}
We iterate this procedure until the difference between $u^k$ and $u^{k+1}$ is sufficiently small.
\end{algorithm}

We continue to describe some heuristic reasons  as follows: 
Note that $T(\bfx) \approx \nabla u^{k+1}(\bfx)$.
In particular, when $\bfx= v_\theta\in 
\partial V$, $\nabla u^{k+1}(v_\theta)= w_\theta\in \partial W$ by the boundary
condition in (\ref{newboundarycondition}). Our choice of 
$w_\theta$ associated with $v_\theta$  will make the distance between each other the 
smallest when the three points $x_0, \bfx, \nabla u^{k+1}(\bfx)$ are on the same
line which contains the points $v_\theta$ as illustrated in Figure~\ref{ex1}.
Furthermore, if $u^{k+1}$ is 
a convex function, then $D_x u^{k+1}$ is increasing and so is $D_y u^{k+1}$. 
In this case, when $\bfx$ inside $V$ starting from $x_0$ increases to 
$v_\theta\in \partial V$, 
$\nabla u^{k+1}(\bfx)$ also increases from $y_0$ to $w_\theta\in \partial W$. 
Our conditions (\ref{newboundarycondition}) and (\ref{normalization}) will map $V$ to $W$. 

We will demonstrate that our boundary conditions can build up a transportation boundary map which 
has a good chance to solve the optimal transport problem in a setting where both 
$V$ and $W$ are star-shaped domains, e.g. convex domains with their geometric centers.

In the following, we shall present a few  computational examples to 
demonstrate the performance of our center matching algorithm. 

\subsubsection{Accuracy of Algorithm~\ref{alg2}}
Let us use the known solution of a transportation problem to check the accuracy 
of our Algorithm~\ref{alg2}. More precisely, we aim to solve a transport problem and compare our 
results with those obtained in previous works \cite{BFO14, LR17, HCX22}. 

\begin{example}\label{optimalex4} 
In this example, the problem involves mapping a given density function $f(x_1,x_2)$ over the square $(-0.5,0.5)\times(-0.5,0.5)$ onto a uniform density function $g=1$ within the same square. The density function $f(x_1,x_2)$ and the mapping function $q(z)$ are defined as follows:
\begin{equation}
f(x_1,x_2) = 1+4(q''(x_1)q(x_2)+q(x_1)q''(x_2))
+ 16(q(x_1)q(x_2)q''(x_1)q''(x_2)-q'(x_1)^2q'(x_2)^2),
\end{equation}
where
\begin{equation}
q(z) = \left(-\frac{1}{8\pi}z^2 + \frac{1}{256\pi^3}+\frac{1}{32\pi}\right)\cos(8\pi z) + \frac{1}{32\pi^2}z\sin(8\pi z).
\end{equation}
The exact solution to this transport problem is given by:
\begin{equation}
u_{x_1}(x_1,x_2) = x_1 + 4q'(x_1)q(x_2), \quad u_{x_2}(x_1,x_2) = x_2 + 4q(x_1)q'(x_2),
\end{equation}
as explained in \cite{BFO14}. 
To find the solution, we utilize our proposed Algorithm~\ref{alg2} to find 
the resulting spline solution $u_s$ which is accurate, as demonstrated in Figure~\ref{fig:ex3fg}. Indeed, we evaluate the solution at $1000\times 1000$ equally-spaced points over $V$ to compute
the maximum of the error and RMS error against the exact solution as well as the maximal error of the \MAE.
\begin{table}[htpb]
\caption{Numerical results for Example~\ref{optimalex4} using Algorithm~\ref{alg2} based on splines of various degrees \label{table:optimalex4}}
\centering
\begin{tabular}{c c c c c}
\hline
Degree & CPU time & Max error & L$^2$ error & $\|\text{det}(D^2u_k)-f/g(\nabla u_k)\|_\infty$ \\
\hline
5 & 1.75 & 3.49e-03 & 4.48e-04 & 1.12e-01 \\
8 & 3.21 & 3.37e-05 & 3.35e-06 & 1.17e-03 \\
11 & 38.0 & 3.37e-07 & 3.11e-08 & 1.62e-05 \\
14 & 95.9 & 2.96e-09 & 4.87e-10 & 2.98e-08 \\
\hline
\end{tabular}
\end{table}

\begin{figure} 
    \centering
\includegraphics[height=4cm ]{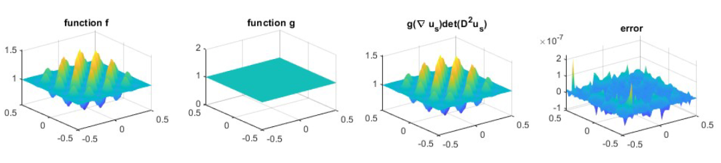}
    \caption{Function $f, g$, our spline solution $f_s$ and error plot $g(\nabla u_s)\text{det}(D^2u_s)-f$ 
    based on Algorithm~\ref{alg2} for Example~\ref{optimalex4} } 
    \label{fig:ex3fg}
\end{figure}

We present detailed computational results in Table~\ref{table:optimalex4}, which summarizes the numerical outcomes obtained using Algorithm~\ref{alg2}. The table includes the degree of the spline functions, the CPU time, the maximum error,  
the RMS error for $L^2$ error, and $|\text{det}(D^2u_k)-f/g(\nabla u_k)|_\infty$ for each degree. 
These results provide valuable insights into the accuracy and efficiency of our method. Both the maximal error and $L^2$ error are
much smaller than the numerical results from the computational method discussed in \cite{BFO14, LR17, HCX22}.

In addition, let us show the exact functions $f$ and $g$ and then $g(\nabla u_s) \det(D^2 u_s)$ and its error
$f-g(\nabla u_s) \det(D^2 u_s)$ in  Figure~\ref{fig:ex3fg}. 
We can see that the error $f-g(\nabla u_s) \det(D^2 u_s)$ is near 0.  

Finally, we mention that the researchers in a recent publication \cite{HCX22} studied this example with various density functions. Our spline collocation method can certainly 
reproduce the results in \cite{HCX22}. The detail is left to the interested reader. 

\end{example}

Certainly, one possible question is how to choose a good center $x_0$. 
It turns out that for this example above, a center $x_0$ makes no significant difference.  
We used various points as a center for our Algorithm~\ref{alg2}.  Numerical results are very similar.  


\begin{example}\label{optimalex5}
In \cite{BFO14}, the researchers explained another example to show their computational
method. They consider $V=W$ which is  the unit square $[-1,1]\times[-1,1]$ 
and two Gaussian density functions. 
We consider two cases. 

Case 1.
The target density is simply a Gaussian in the center of the domain:
	$$ g(y) =  2+\frac{1}{0.2^2}\exp\left(-\frac{0.5|x|^2}{0.2^2}\right).$$
For the source density, we use four Gaussians centered at the four corners of the domain
to form the density function.  For example, in the quadrant $[-1,0]\times[-1,0]$, we use
$$ f(x) =  2+\frac{1}{0.2^2}\exp\left(-\frac{0.5|x-(-1,-1)|^2}{0.2^2}\right).$$
Figure \ref{fig:ex4fg} (the two graphs most left) depicts the density functions $f(x)$ and $g(y)$. In the two graphs most right of Figure \ref{fig:ex4fg}, we present the function $g(\nabla u_s)\det(D^2 u_s)$ obtained using a spline solution $u_s$, along with the plot of the error $g(\nabla u_s)\det(D^2 u_s) - f$. 
\begin{figure}
    \centering
    \includegraphics[height=4cm ]{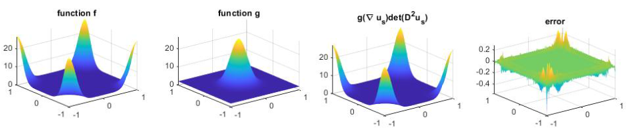}
    \caption{Functions $f, g$, our numerical solution, and error plot $\|g(\nabla u_s)\text{det}(D^2u_s)-f\|_\infty$ for the case 1 in Example \ref{optimalex5}}
    \label{fig:ex4fg}
\end{figure}

Case 2.
In this case, we interchanged the source and target densities from Case 1. The source density $f(x)$ in Case 2 is now the same as the target density in Case 1, and the target density $g(y)$ in Case 2 is the same as the source density in Case 1.
The density functions $f(x)$ and $g(y)$ are depicted in Figure \ref{fig:ex4fg2} (the top two graphs). In the bottom graph of Figure \ref{fig:ex4fg2}, we show the function $g(\nabla u_s)\det(D^2 u_s)$ obtained using a spline solution $u_s$, along with the plot of the error $g(\nabla u_s)\det(D^2 u_s) - f$. 

\begin{table}[htp]
\centering
\begin{tabular}{ccccc}
\hline
Algorithm &$(x_0, y_0)$ &  Time(s)& $|\det(D^2u_s) - \frac{f}{g(\nabla u_s)}|_2$ & $I[\nabla u_s]$ \\ 
\hline
Algorithm~\ref{alg2} & (0.0,0.0)& 163.05 & 6.07e-03 & 1.3115550  \\   
\hline
\end{tabular}
\caption{Numerical results for Example \ref{optimalex5} using spline functions of $D=12, r=2$, where $I[\nabla u_s]=\int_V |x-\nabla u_s(x)|^2f(x) dx$. }\label{table:example4}
\end{table}

\begin{figure}
    \centering
    \includegraphics[height=3.5cm ]{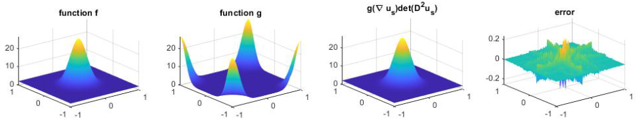}
    \caption{Functions $f, g$, our numerica solution of $f$, and error plot $g(\nabla u_s)\text{det}(D^2u_s)-f$ for the case 2 in Example \ref{optimalex5}}
    \label{fig:ex4fg2}
\end{figure}

Furthermore, Figure \ref{fig:ex4u} displays the Brenier potential $u$ for both Case 1 and Case 2. 
Additionally, Table \ref{table:example4} provides numerical results obtained using two algorithms with $D=12$ and $r=2$, 
where $I[\nabla u_s]=\int_V |x-\nabla u_s(x)|^2f(x) dx$.

\begin{figure}
  \centering
  \begin{subfigure}[b]{0.25\textwidth}
    \includegraphics[width=\textwidth]{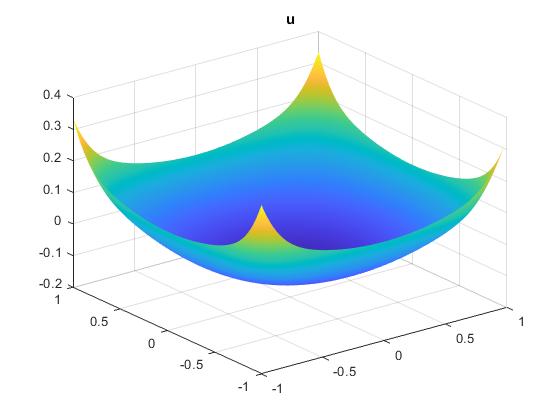}
  \end{subfigure}
  \begin{subfigure}[b]{0.25\textwidth}
    \includegraphics[width=\textwidth]{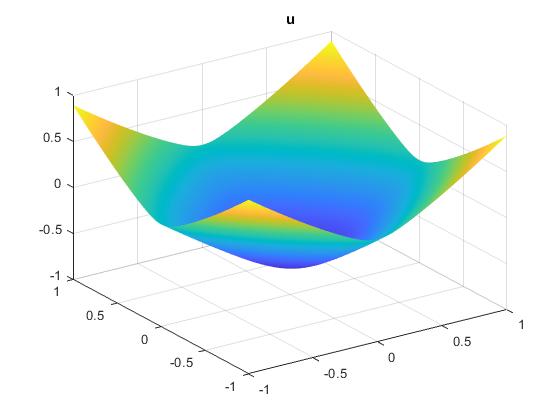}
  \end{subfigure}
  \caption{Brenier potential $u$ for the case 1(left) and case 2(right)}\label{fig:ex4u}
\end{figure}

\end{example}

\subsubsection{A Natural Transportation Problem}  
From the numerical examples in the previous subsections, we can observe that our spline-based collocation method performs very well in obtaining numerical solutions 
for the \MAE numerically. Now, we will utilize this method to solve a special case of the optimal transportation problem. 
In this subsection, we will numerically compute the solution to the \MAE  
with a constant on the right-hand side. This problem is related to a natural transportation problem because $\frac{f(\mathbf{x})}{g(\nabla u(\mathbf{x}))}\equiv 1$ for all $\mathbf{x}\in V$. To find a numerical solution, we will employ Algorithm~\ref{alg2}.

We claim that the $u_s$ is the solution to the optimal transport problem if the solution $u_s$ is a convex function,  the map $\nabla u_s$ transforms $V$ to the entire domain $W$ and $u_s$
satisfies the \MAE  with the right-hand side $f\equiv 1$. Indeed, according to Theorem~\ref{BrenierThm1}, the solution $u$ is unique up to a constant, such that the push-forward $\nabla u$ from $V$ to $W$ satisfies $\nabla u (V)=W$ and the solution of the \MAE $ u$ is convex.

\begin{figure}[thpb]
 \centering
\includegraphics[width = 0.75\textwidth]{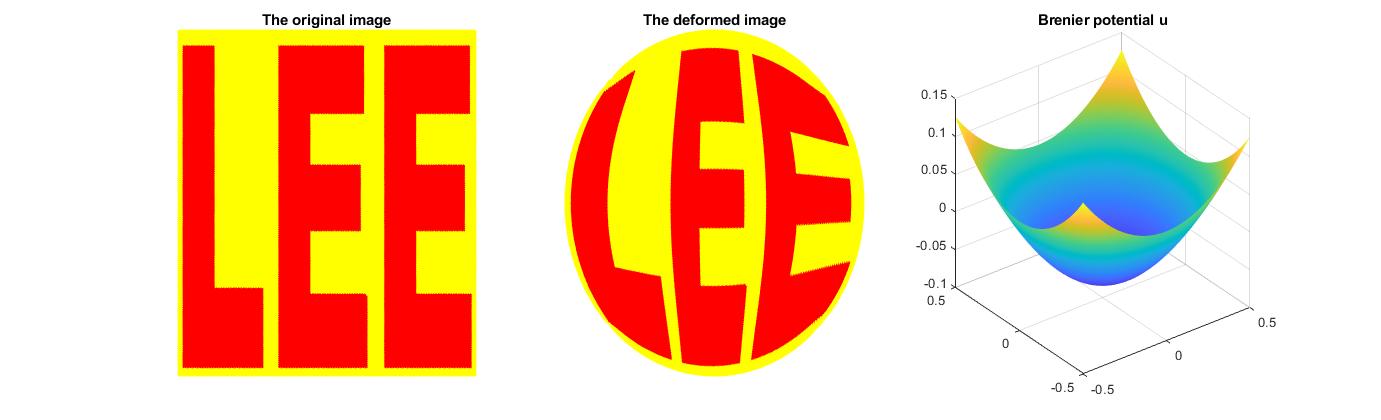} 
\caption{Deformation of an image (left) from a square domain to a circle domain (middle) together with the Brenier potential (right). The moving cost is $I[\nabla u]=5.4490e-03$}
\label{lfig1}
\end{figure}

To demonstrate that $\nabla u(V)=W$, we use the image representing a density function $f$ over a rectangular domain $V$ and plot the pixel values $f(x)$ as $g(\nabla u(x))$ at $\nabla u(x)\in W$. If the image $f$ is entirely transformed to $W$, then it indicates that $\nabla u(V)=W$. Figure~\ref{lfig1} 
illustrate that the image from the rectangular domain $V$ is perfectly mapped 
onto the disk $W$ (in the middle of Figure~\ref{lfig1}, and the resulting solution $u$, known as the Brenier potential, 
is convex as depicted on the right of Figure~\ref{lfig1}. 

\begin{figure}[thpb]
\centering
\includegraphics[width = 0.75\textwidth]{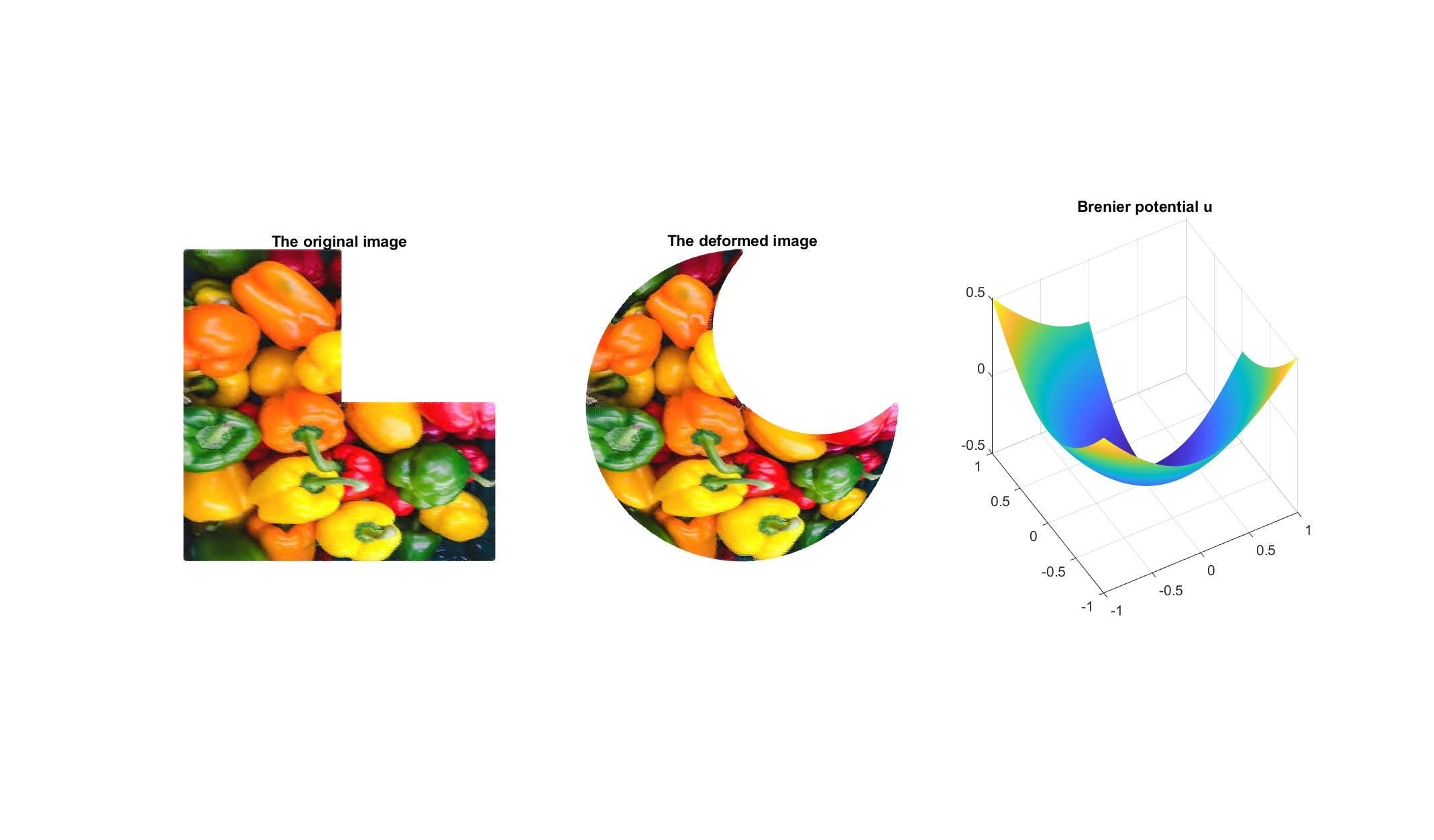} 
\caption{Deformation of an image (left) from an L-shaped domain to a new L-shaped (moon-like) domain (middle) together with the Brenier potential (right) }
\label{lfig2}
\end{figure}

\begin{figure}[thpb]
\centering
\includegraphics[width = 0.8\textwidth]{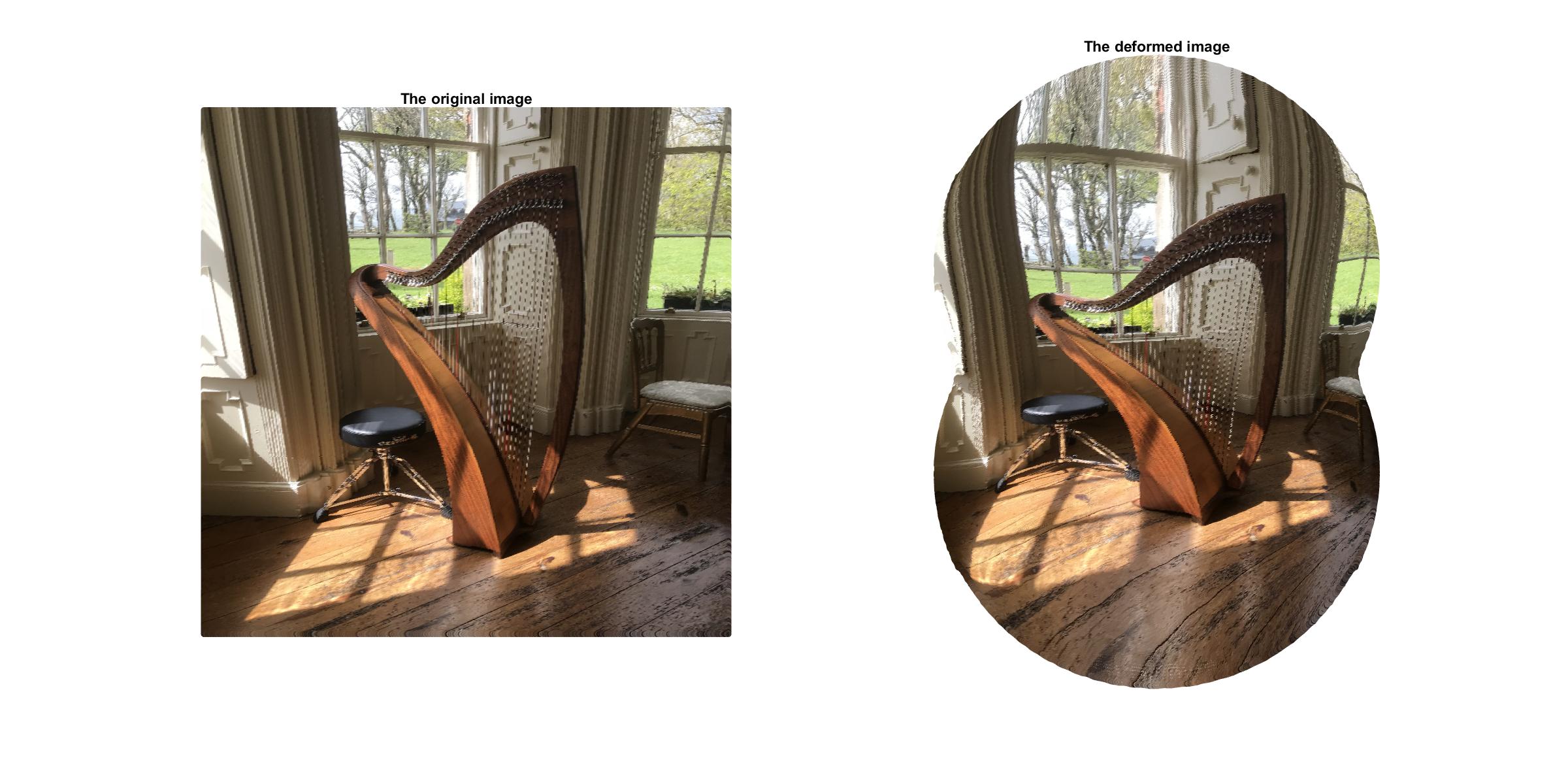} 
\caption{Deformation of an image (left) from a square domain to a nonconvex domain. The Brenier potential 
is convex which is not shown above.}
\label{lfig3}
\end{figure}

The following are some examples from our computational algorithm. See Figures~\ref{lfig1}--\ref{lfig3}.

Similarly, we use Algorithm~\ref{alg2} to find the solution to the natural transportation problem
 to any pair of star-shaped domains that are not 
convex. In Figure~\ref{lfig2}, we use an L-shaped domain and a moon-shaped domain and Algorithm~\ref{alg2} finds $u_s$ which is convex as shown on the right-hand of Figure~\ref{lfig2}. The middle of Figure~\ref{lfig2} 
shows that $\nabla u_s(V)=W$.    
We are able to generate many image transforms from one 
star-shaped domain to another. The details are omitted. The interested reader may contact the authors for a MATLAB code. 

\subsection{Applications of the Optimal Transport Problem}
With the computational algorithm, i.e. Algorithm~\ref{MAEalg} with the center matching algorithm~\ref{alg2}, 
we are now able to solve some optimal transport problems. 
A special optimal transportation problem is to 
to solve the Monge-Amp\'ere equation
\begin{equation}
\label{MAE2b}
\left\{ \begin{array}{clc}
\det (D^2 u)(\bfx) &=f(\bfx)/A(f), & \quad \bfx \in V,\cr
\nabla u(V) &= V, & u \hbox{ is convex over } V, 
\end{array}
 \right.  
\end{equation}
where $f(\bfx)\ge f_0>0$ and $g(\bfy)\equiv D(f)/A(W)$ are  given with 
$D(f)$ being the volume of density $f$ over $V$ and $A(W)$ is the area of $W$.  Note that we can 
approximate $A(f)$ accurately by using bivariate splines easily. Three applications will be explained below.  

\begin{figure}[htpb]
\includegraphics[width = 1\textwidth, height=0.5\textwidth]{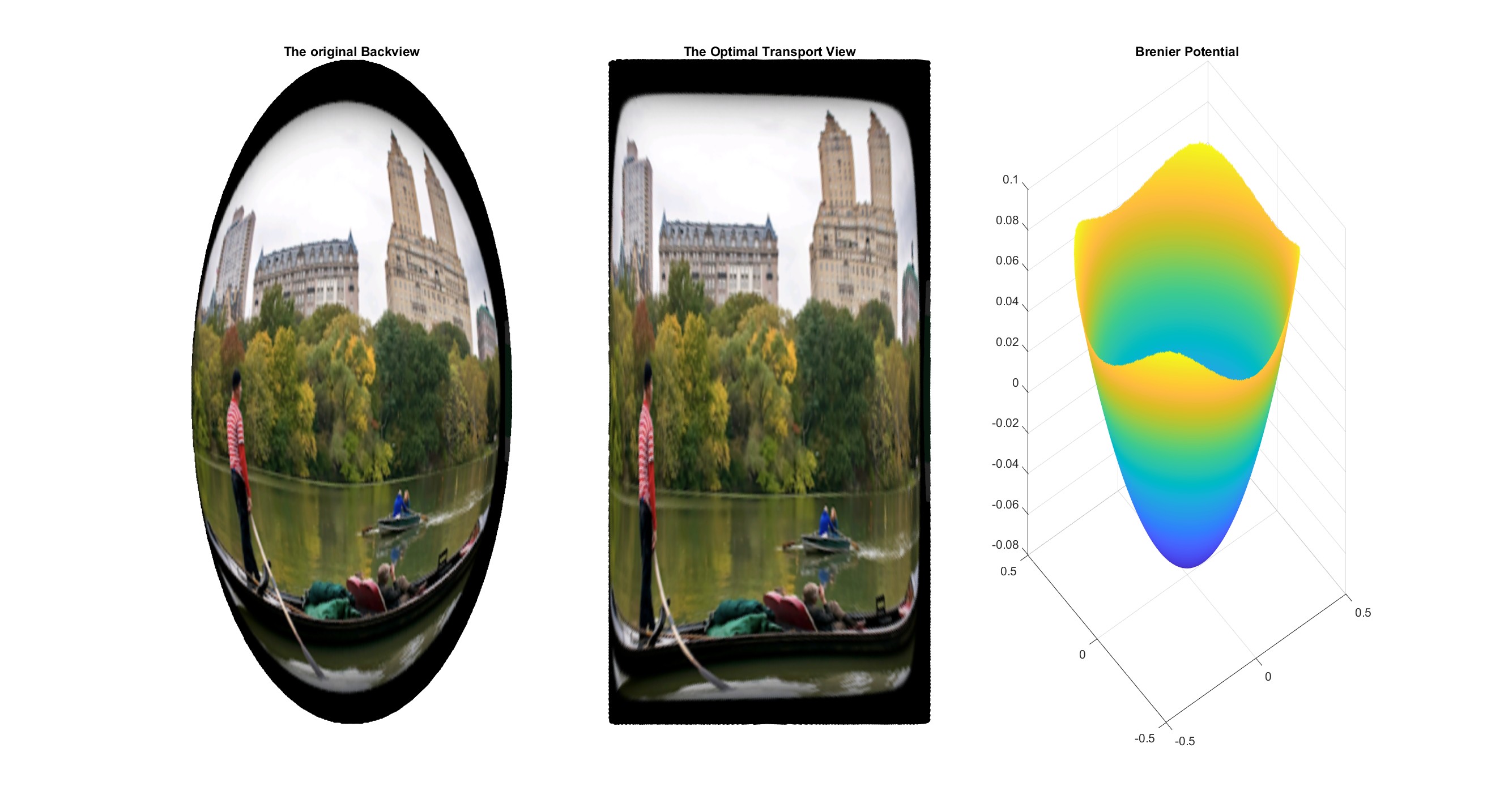} 
\caption{A fisheye view image (left), the optimal transport image (middle), and
the Brenier potential which is convex. The RMSE of the \MAE is equal to 
$1.124499e-3$. \label{LaisEx4}} 
\end{figure}

Let us first discuss how to transport a fisheye image from a surveillance camera onto a standard rectangular image. 
Clearly, a fisheye image is a density function over a circular domain $V$ centered at the origin with 
radius $1$. We shall choose a constant function 
$f$ and $g=pi/4$ which is the volume of $f$ over $D$. Assume $W$ is a square domain $[-0.5, 0.5]^2$. We use our method to find the Brenier potential 
$u$. Then $\nabla u$ is the onto map as shown in Figure~\ref{LaisEx4}. For any pixel value $z(x,y)$ with $(x,y)\in V$, we let $z(x,y)$ be the value 
associated with $\nabla u(x,y)$ on $W$.  
In this way, we map the image on $V$ to an image over the rectangular domain $W$.   Based on the images, we 
are able to show the bijectivity of the map   $\nabla u$. To show the convexity of
Brenier potential, we plot the function $u$. In addition, we also present the root mean square error (RMSE) of the \MAE over all image pixel locations $(x,y)$ 
\begin{equation}
\label{MAE2}
E(x,y)= \det( D^2 u)(x,y) - f(x,y)/g(\nabla u(x,y)
\end{equation}
to show the accuracy of our spline solution.  Note that the size of the pixel value locations is more than $501
\times 501$.  

The Brenier theorem mentioned at the beginning of the paper tells us that the
spline solution is a very good approximation of the optimal transport problem. We can see that the man standing straight 
on the boat which is horizontal, and the river bank is horizontal. Four buildings are straight except for the two tip parts of
the building on the top-right corner of the image which are bent. More research is needed to make them straight upward. 

\begin{figure}[htpb]
\centering
\includegraphics[width = 1\textwidth, height=0.5\textwidth]{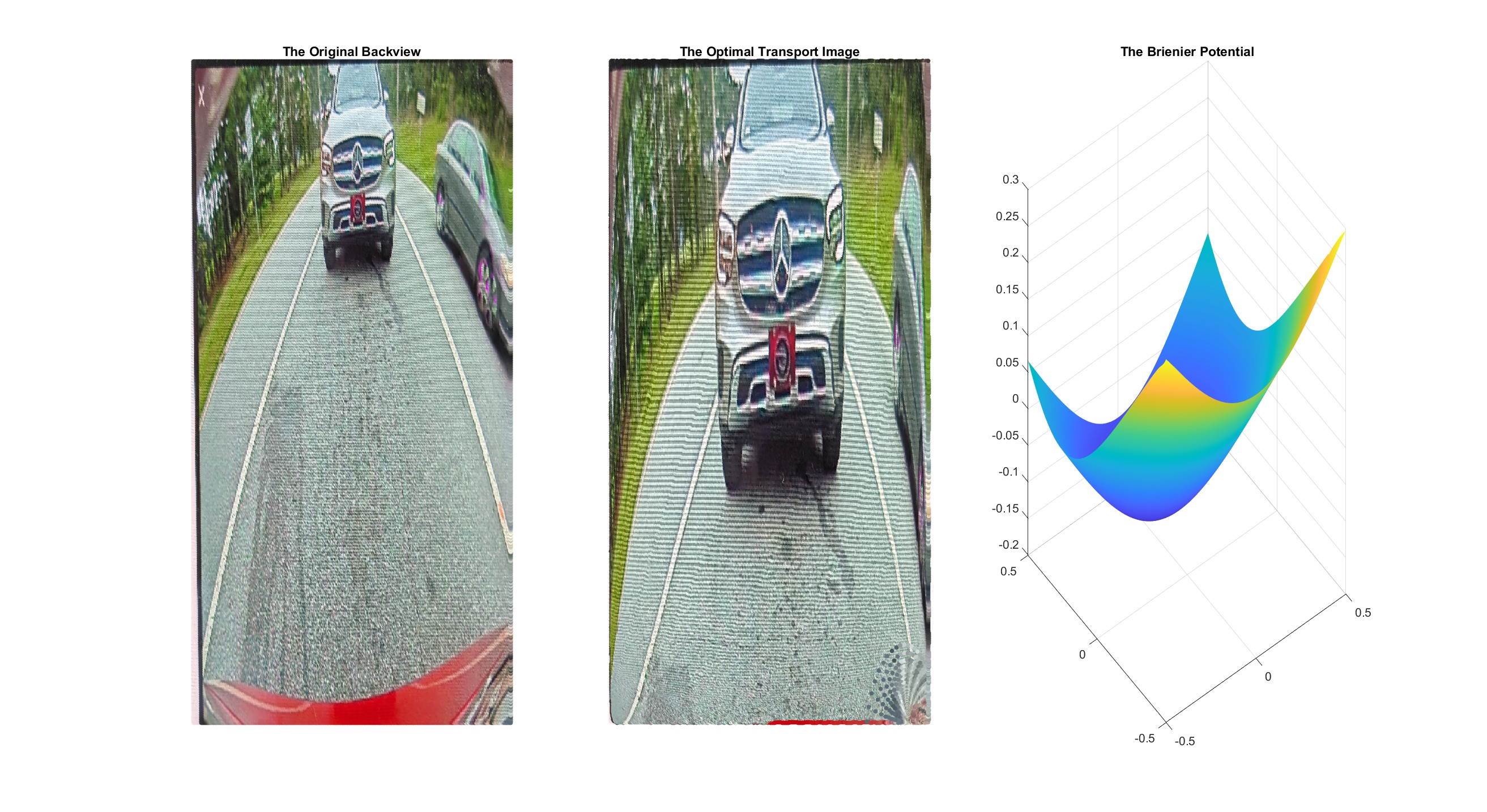} 
\caption{A standard back view image (left), the optimal transport image (middle), and
the Brenier potential which is convex. The RMSE of the \MAE is equal to 
$6.037135971841e-3$. \label{LaisEx5}} 
\end{figure} 

Our second example is to deform the back view image of a car. Note that 
$V=W$. 
We used the following function 
\begin{equation}
\label{density}
f(x,y)=  \exp(\alpha (x-t)^2+\beta (y-s)^2)
\end{equation}
with appropriate parameters $\alpha, \beta, t, s$.  We choose $g$ to be 
a constant. The top central part of
the image needs an enlargement and hence $s\not=0$.  
Again we use the Brenier theorem mentioned at the beginning of the paper to see that the
spline solution is a very good approximation of the optimal transport problem. The car in the optimal transport 
image is bigger and looks closer to the real-life situation as the distance between the two cars is indeed very close.

\begin{figure}[htpb]
\centering
\includegraphics[width = 1\textwidth, height=0.5\textwidth]{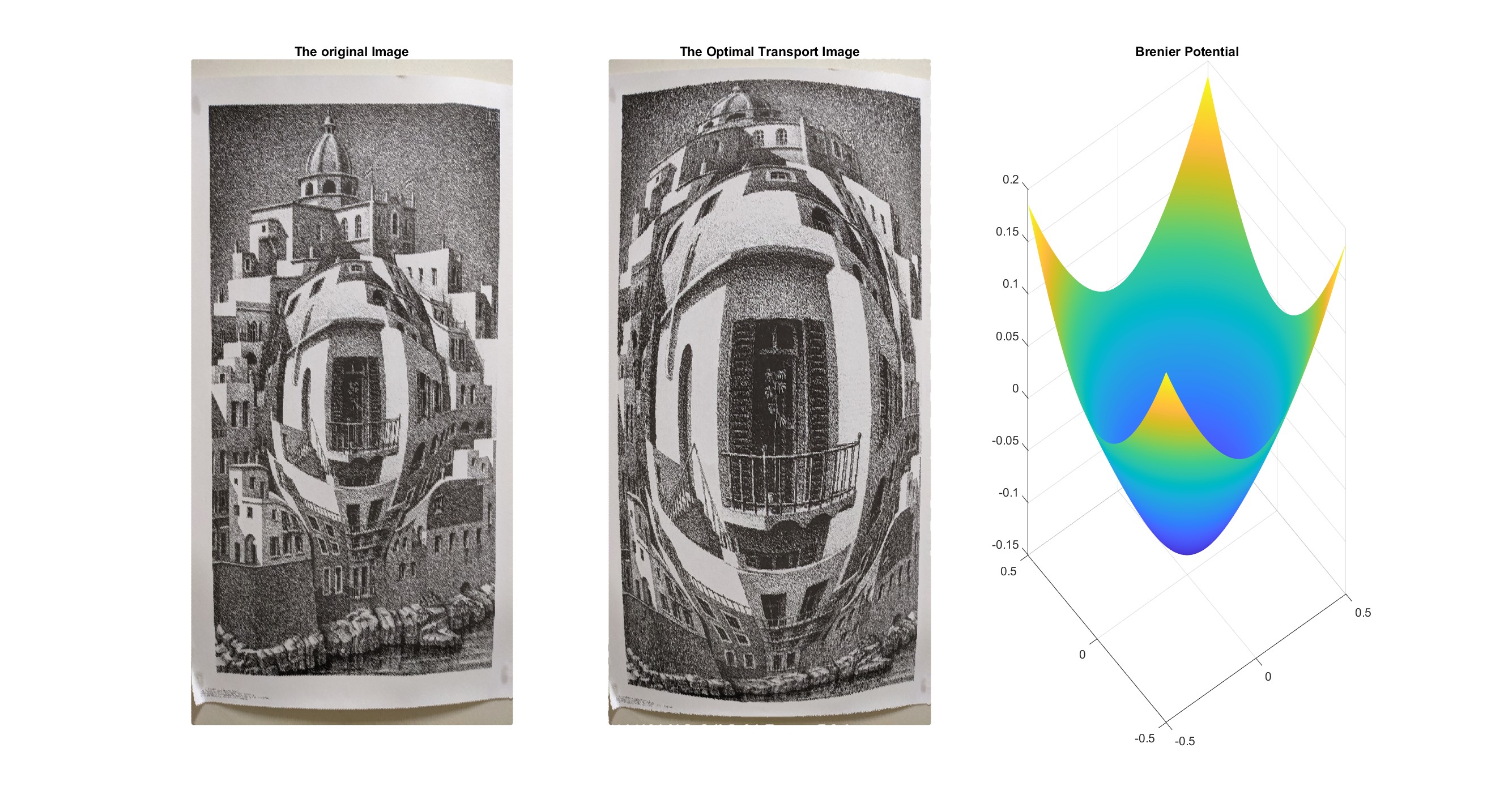} 
\caption{An Escher art (left), the optimal transport image (middle), and
the Brenier potential which is convex. The RMSE of the \MAE is equal to 
$1.351757e-3$.}\label{LaisEx6} 
\end{figure}

The experimental study reminds the first author that M.C. Escher created an 
image as shown on the left of Figure~\ref{LaisEx6}.  We can do the same thing
for any image to have a deformed view at the center of the image. Our spline method can generate such a view easily for 
any image. See Figure~\ref{LaisEx6} for the Brenier potential which is convex and the RMSE
of the \MAE in the caption of the figure.       
Note that the map $\nabla u(x,y)$ can be used to generate an Escher art-like image for any given image. We can also 
adjust these parameters $\alpha, \beta, t, s$ to have many different images similar to the Escher art.

\begin{figure}[thpb]
\centering
\includegraphics[width = 1\textwidth, height=0.5\textwidth]{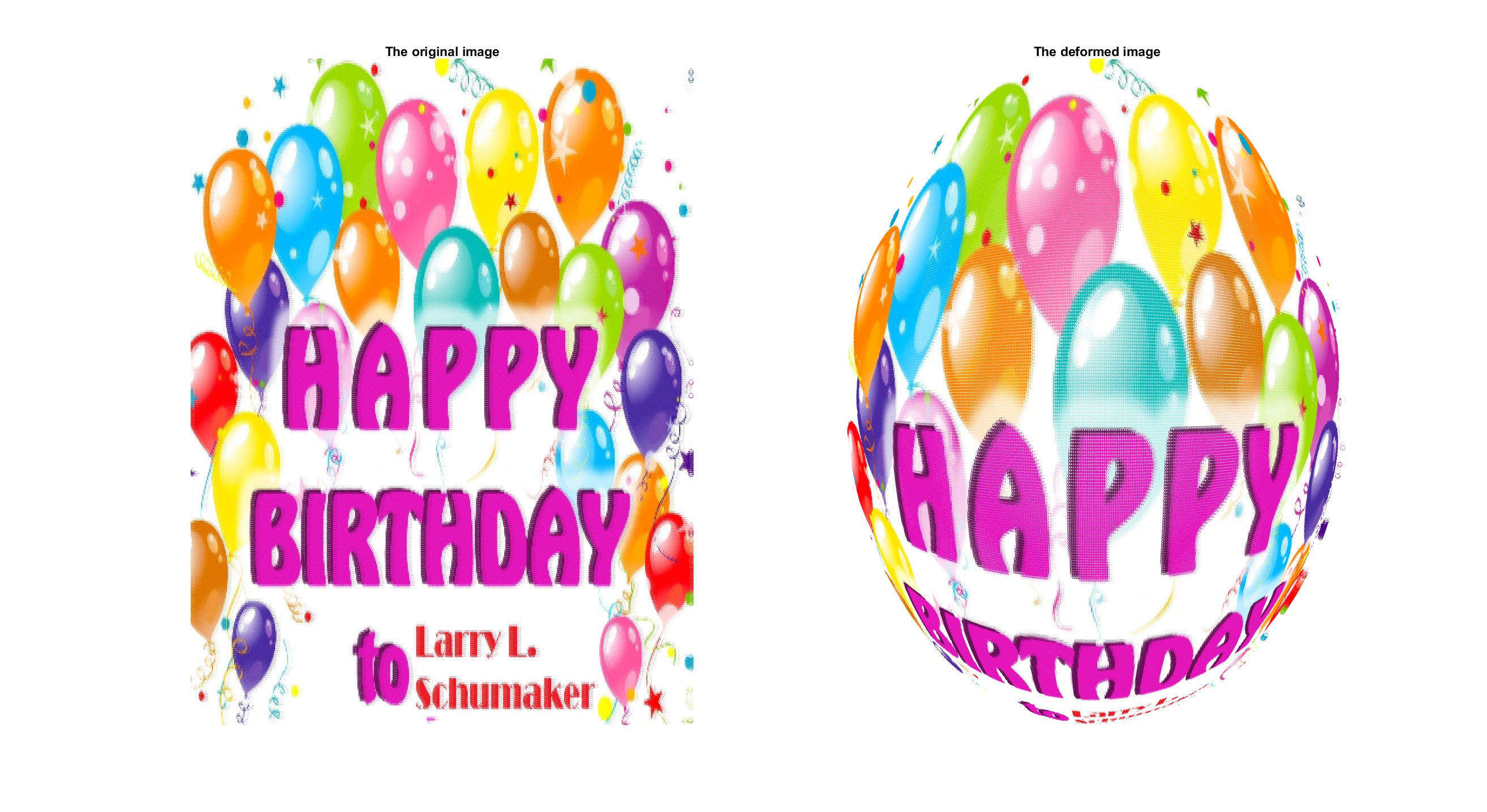} 
\caption{Deformation of an image (left) to the image on the right which looks like a 3D view of the balloon. }
\label{LaisEx1}
\end{figure}

\begin{figure}[thpb]
\centering
\includegraphics[width = 0.5\textwidth]{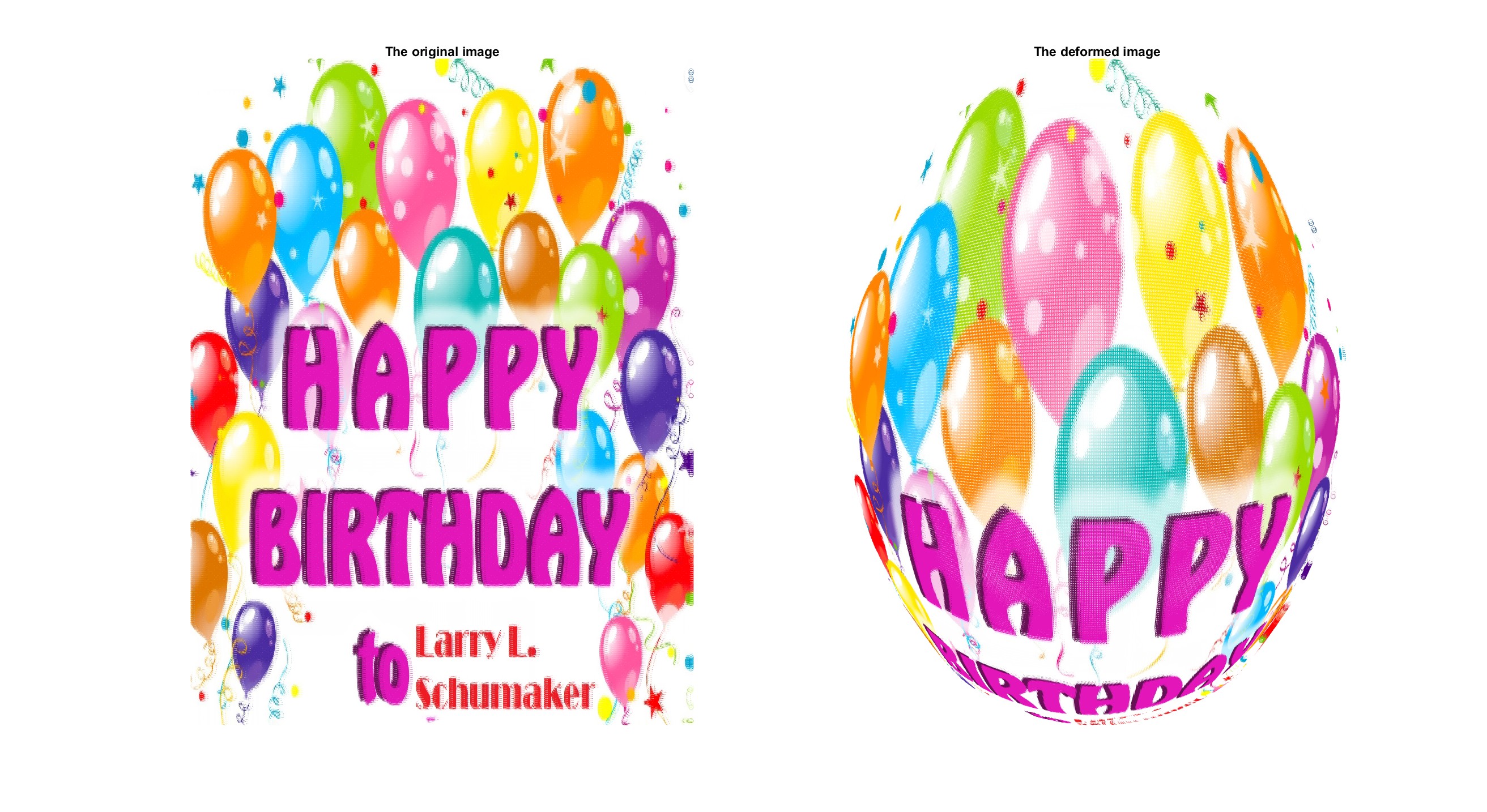} 
\includegraphics[width = 0.5\textwidth]{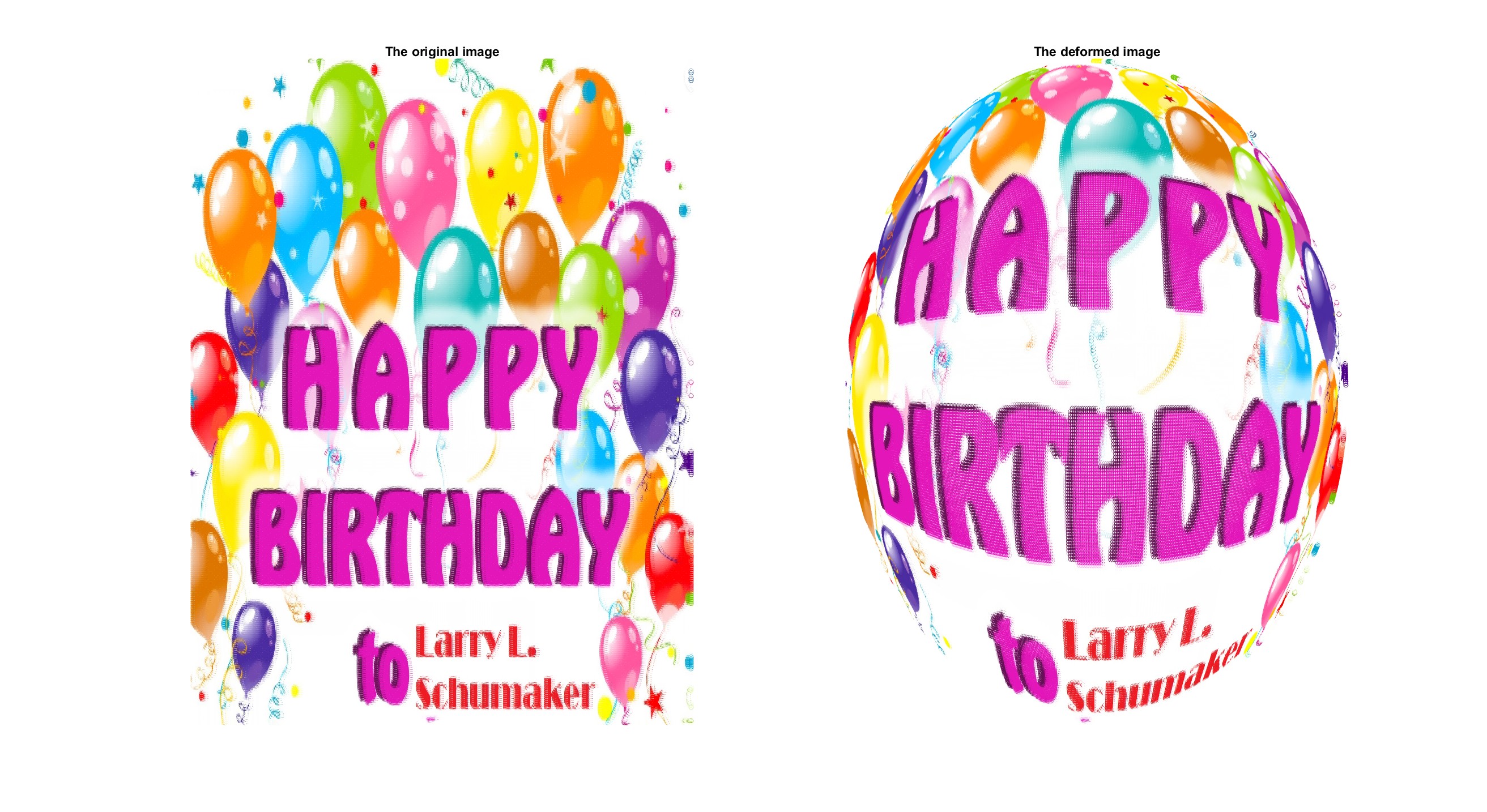} 
\caption{Deformation of an image (left) to the image on the right which looks like a 3D view of the balloon. }
\label{LaisEx2}
\end{figure}

Finally, we present another example. 
Let $V$ be a square domain and $W$ be an oval domain. For various density functions $f$ in 
(\ref{density}) on $V$ and the constant function $g$ on $W$ with $g= \int_V f(x,y)dxdy$, we solve the  \MAE \ref{MAE2b}. 
Note that we use bivariate splines to find a spline approximation of $f$ so that we can find 
$g= \int_V f(x,y)dxdy$ easily based on the integration formula in \cite{LS07}.  
Our computational algorithm produces a  deformation of
image $B$ over $V$ to an image
$D$ over $W$,  It is a bijective map with a convex Briener function which is not shown.  
By using various parameters 
$t$ in (\ref{density}), we generate more deformed images as shown in Figures~\ref{LaisEx1} and 
\ref{LaisEx2}.

\section{Appendix}
We explain how to solve the Poisson equation by using bivariate spline collocation method. Recall   
\begin{equation}
\label{Poisson}
-\Delta u  = f, \quad {\hbox{ in }} \Omega\subset \mathbb{R}^d, \hbox{ and }  
   u = g,  \quad \hbox{ on } \partial \Omega.
\end{equation}
When $\Omega$ has a uniform positive reach, the solution to the Poisson equation will be 
in $H^2(\omega)$ according to \cite{GL20}. 
In \cite{LL22}, we have developed  the multivariate spline-based collocation method to solve linear partial 
differential equations including the Poisson equation above.  We mainly use $C^r$ spline functions with $r\ge 1$ to approximate the solution $u$ using the PDE directly without any weak formulation. 

For our method, we shall use the so-called domain points (cf. \cite{LS07} and \cite{LL22}) 
to be the collocation points.   Letting  $\xi_i, i=1, \cdots, N$ 
be the domain points of $\triangle$ and degree $D'>0$, 
where $D'$ will be different from $D$, 
our multivariate spline-based collocation method is to seek a spline function $s
\in S^r_D(\triangle)$ satisfying 
\begin{equation}
\label{Poisson2}
-\Delta s(\xi_i)  = f(\xi_i), \quad \xi_i\in \Omega, i=1, \cdots, N, \hbox{ and }  
   s(\xi_j)  = g(\xi_j),  \quad \xi_j\in \partial \Omega, j=1, \cdots, M.
\end{equation}
Instead of using locally supported basis functions for 
spline space $S^r_D(\triangle)$, we will begin with a discontinuous  spline space  $s\in S^{-1}_D(\triangle)$
and then add the smoothness conditions.  That is, we write $s= \sum_{i=1}^N c_i \phi_i$ with discontinuous polynomial 
functions $\phi_i$ as basis functions. When $s\in S^r_D(\triangle)$, the coefficients $c_i$ satisfy the so-called 
smoothness conditions.  
Note that the smoothness conditions across each interior edge of $\triangle$ 
can be written as a system of linear equations, i.e. $H\bfc =0$, where $\bfc$  
is the coefficient vector of spline function $s\in S^{-1}_D(\triangle)$ such that a spline function $s$ with coefficient vector
$\bfc$ will be in $C^r(\Omega)$.  We refer to \cite{LS07} for detail and \cite{ALW06} for implementation.  
One of the key ideas is to let a computer decide how
to choose ${\bf c}$ to satisfy $H{\bf c}=0$ and (\ref{Poisson2}) above simultaneously. 
With  $s= \sum_{i=1}^N c_i \phi_i$, the first equation in (\ref{Poisson2}) can be written as a matrix format, say $-K \bfc
=\bff$ with $\bff$ being the vector on the right-hand side of the equation. 
 Our collocation method for (\ref{Poisson}) is to find ${\bf c}^*$ by solving the following constrained minimization: 
 \begin{align}
\label{min1}
\min_{\bf c} J(c)=\frac{1}{2}(\alpha \|B{\bf c}- {\bf g}\|^2+\beta 
\|H{\bf c}\|^2)  
\quad \text{subject to } \|K{\bf c} + {\bf f}\| \le \epsilon_1,
\end{align}
where $B, {\bf g}$ are associated with the boundary condition, $H_r$ is associated with 
the smoothness condition with $r=2$, $\alpha>0, \beta>0$ are fixed parameters, and $\epsilon_1>0$
is a given tolerance. 

To explain the convergence of the numerical solution $s^*$ with coefficient vector $\bfc^*$, we recall the following 
fundamental result.
\begin{lemma}([Lai and Schumaker, 2007\cite{LS07}])
\label{lem1}
Let $k \geq 3r+2$ with $r\geq 1$. 
Suppose $\triangle$ is a quasi-uniform triangulation of $\Omega$. 
Then for every $u\in W_q^{k+1}(\Omega),$ there exists a quasi-interpolatory spline $s_u\in \mathcal{S}^r_k (\triangle)$ such that 
\begin{eqnarray*}
	\|D^\alpha_x D^\beta_y (u- s_u)\|_{q, \Omega}\le C |\triangle|^{k+1-\alpha-\beta} |u|_{k+1,q,\Omega}
\end{eqnarray*}
for a positive constant $C$ dependent on $u, r, k$ and the smallest angle of $\triangle $, and for all $0\leq 
\alpha+\beta \leq k$ with 
\begin{eqnarray*}
	|u|_{k,q,\Omega}:=(\sum_{a+b=k}\|D_x^a D_y^b u\|^q_{L^q(\Omega)})^{\frac{1}{q}}.
\end{eqnarray*}
\end{lemma}

Suppose that the solution $u$ to the Poisson equation is in $H^3(\Omega)$. According to the result above, there exists a 
quasi-interpolatory spline $s_u$ such that 
$$
\|\Delta u- \Delta s_u\|_{L^2(\Omega)} \le \epsilon 
$$
if the size $|\triangle|$ of the underlying triangulation $\triangle$ of $\Omega$ is mall enough.  Also, if we choose $D'$ 
large enough such that the root mean square error $\|\Delta u- \Delta s_u\|$ based on the collocation points is also 
less than or equal to $\epsilon_1\ge \epsilon$. Now the feasible set of (\ref{min1}) is not empty. 
Since the minimization (\ref{min1}) is a convex
minimization problem over a convex nonempty feasible set, the problem (\ref{min1}) will have a
unique solution. We have developed a computational algorithm to solve the minimization (\ref{min1}). 

Next, we need to solve that the spline $s^*$ with the minimizer coefficient vector $\bfc^*$ of (\ref{min1}) is a good 
approximation of the solution $u$.   
That is, recall we often use the standard $H^2$ norm 
\begin{equation}
\label{H2norm}
\|u\|_{H^2} =\|u\|_{L^2(\Omega)}+ 
\|\nabla u\|_{L^2(\Omega)}+ \sum_{i,j=1}^d 
\|\frac{\partial}{\partial x_i}\frac{\partial}{\partial x_j}u\|_{L^2(\Omega)}
\end{equation}
for all $u$ on $H^2(\Omega)$.
In this paper and in \cite{LL22}, we let
$$\|u\|_L= \|\Delta u\|_{L^2(\Omega)}, ~\forall u \in H^2(\Omega)\cap H^1_0(\Omega)$$
be a norm on $H^2(\Omega)\cap H^1_0(\Omega).$ Then 

\begin{lemma}
\label{normequiv}
 There exists two positive constants $A>0$ and $B>0$ such that  
\begin{equation}
\label{normequivalence}
A\|u\|_{H^2} \le  \|u\|_L \le B \|u\|_{H^2}.
\end{equation}
\end{lemma}

See proof in \cite{LL22}. Suppose that $u- s^* \equiv 0$ on $\partial \Omega$, we have 
$$
\|u - s^*\|_{H^2(\Omega)} \le \frac{1}{A} \|u-s^*\|_{L^2(\Omega)} = \frac{1}{A} \| \Delta  u - \Delta s^*\| \le 
\frac{\epsilon_1}{A}.
$$
This explains that the minimization (\ref{min1}) can be used to find a good approximation of the 
solution $u$ to the Poisson equation.  
\end{document}